\documentclass{article}
\usepackage[english]{babel}
\usepackage[utf8]{inputenc}
\usepackage[T1]{fontenc}
\usepackage{amsfonts}
\usepackage{amsmath}
\usepackage{amsthm}
\usepackage{amssymb}
\usepackage{hyperref}
\usepackage{enumitem}
\usepackage{array}
\usepackage{mathtools}
\usepackage{thmtools, thm-restate}
\usepackage{xpatch}
\usepackage{tocloft}

\newtheorem{theoremIntro}{Theorem}
\newtheorem{theorem}{Theorem}[section]
\newtheorem*{theorem*}{Theorem}

\newtheorem{corollary}[theorem]{Corollary}
\newtheorem{corollaryIntro}[theoremIntro]{Corollary}
\newtheorem{proposition}[theorem]{Proposition}
\newtheorem{propositionIntro}[theoremIntro]{Proposition}

\newtheorem*{conjecture*}{Conjecture}
\newtheorem{lemma}[theorem]{Lemma}
\theoremstyle{definition}
\newtheorem*{definition}{Definition}
\theoremstyle{remark}
\newtheorem*{remark}{Remark}
\newtheorem*{example}{Example}
\newtheorem*{examples}{Examples}

\newcommand{\C}{\mathbb{C}}
\newcommand{\R}{\mathbb{R}}
\newcommand{\N}{\mathbb{N}}
\newcommand{\Z}{\mathbb{Z}}

\newcommand{\HNC}{\mathbb{H}^{n}_{\C}}
\newcommand{\PHNC}{\partial\mathbb{H}^{n}_{\C}}
\newcommand{\HNR}{\mathbb{H}^{n}_{\R}}

\newcommand{\PU}{\mathrm{PU}}
\newcommand{\PO}{\mathrm{PO}}
\newcommand{\PNC}{\C \mathbb{P}^{n}}

\DeclareMathOperator{\Id}{\text{Id}}
\DeclareMathOperator{\Ricci}{\text{Ricci}}

\newcommand{\tr}{\text{tr}}

\title{Holomorphic functions on geometrically finite quotients of the ball}

\author{William Sarem}
\date{December 27, 2025}

\begin{document} \maketitle 
	
	\begin{abstract}Let $\Gamma$ be a discrete and torsion-free subgroup of $\PU(n,1)$, the group of biholomorphisms of the unit ball in $\C^{n}$, denoted by $\HNC$. We show that if $\Gamma$ is Abelian, then $\HNC/\Gamma$ is a Stein manifold. If the critical exponent $\delta(\Gamma)$ of $\Gamma$ is less than 2, a conjecture of Dey and Kapovich predicts that the quotient $\HNC/\Gamma$ is Stein. We confirm this conjecture in the case where $\Gamma$ is parabolic or geometrically finite. We also study the case of quotients with $\delta(\Gamma)=2$ that contain compact complex curves and confirm another conjecture of Dey and Kapovich. We finally show that $\HNC/\Gamma$ is Stein when $\Gamma$ is a parabolic or geometrically finite group preserving a totally real and totally geodesic submanifold of $\HNC$, without any hypothesis on the critical exponent.
	\end{abstract}

	In this article we study the existence of non-constant holomorphic functions on quotients of the complex hyperbolic space $\HNC$ of dimension $n$, thought of as the unit ball in $\C^{n}$, by discrete and torsion-free subgroups $\Gamma$ of $\PU(n,1)$. More precisely, we give sufficient conditions for $\HNC/\Gamma$ to be holomorphically convex or Stein. These conditions involve the group structure of $\Gamma$ or its \textit{critical exponent} $\delta(\Gamma)$, which is defined by\begin{equation}\label{eq-def-delta}
		\delta(\Gamma) := \inf \{s \in \R_{+} \ \vert \ \sum_{\gamma \in \Gamma}e^{-sd(o,\gamma o)}< \infty\},
	\end{equation}
	where $o$ is an arbitrary point of $\HNC$ and $d$ the complex hyperbolic distance on the ball, normalized so that the associated Riemannian metric has sectional curvature pinched between $-4$ and $-1$. This number $\delta(\Gamma)$, which does not depend on the choice of $o\in \HNC$, has been first related with the analytical properties of $\HNC/\Gamma$ by Dey and Kapovich in \cite{deyNoteComplexhyperbolicKleinian2020}. These authors have proposed the following conjecture, and have solved it for convex-cocompact subgroups of $\PU(n,1)$.
	
	\begin{conjecture*}[{\cite{deyNoteComplexhyperbolicKleinian2020}}]
		Let $\Gamma$ be a discrete and torsion-free subgroup of $\PU(n,1)$. If $\delta(\Gamma)<2$, then $\HNC/\Gamma$ is a Stein manifold.
	\end{conjecture*}
	
	Our first main result confirms in particular this conjecture for \textit{geometrically finite }subgroups of $\PU(n,1)$.
	
	\begin{theoremIntro}\label{thm-geom-fini-delta-stein}
		Let $\Gamma$ be a geometrically finite and torsion-free subgroup of $\PU(n,1)$.\begin{enumerate}[label=(\alph*)]
			\item If $\Gamma$ is Gromov-hyperbolic, then $\HNC/\Gamma$ is holomorphically convex.
			\item If $\delta(\Gamma)<2$ or if $\Gamma$ preserves a totally real and totally geodesic submanifold of $\HNC$, then $\HNC/\Gamma$ is a Stein manifold.
		\end{enumerate}
	\end{theoremIntro}

	This result applies in particular to free groups, and we obtain the following corollary.
	\begin{corollaryIntro}\label{cor-free}
		Let $\Gamma$ be a geometrically finite free group of $\PU(n,1)$. Then $\HNC/\Gamma$ is a Stein manifold.
	\end{corollaryIntro}
	
	We recall that a parabolic subgroup of $\PU(n,1)$ is a subgroup of $\PU(n,1)$ that fixes a point at infinity and does not contain any hyperbolic element. To establish Theorem \ref{thm-geom-fini-delta-stein}, we need to understand when the quotient
	of the complex hyperbolic space by a discrete parabolic subgroup is Stein. For unipotent parabolic subgroups this has been done in \cite{miebachQuotientsBoundedHomogeneous2024}, and we settle the general case, obtaining the following result.
	\begin{theoremIntro}\label{thm-parabolique-abelien-delta-stein}
		Let $\Gamma$ be a discrete and torsion-free parabolic subgroup of $\PU(n,1)$. \begin{enumerate}[label=(\alph*)]
			\item If $\delta(\Gamma)<2$ or if $\Gamma$ preserves a totally real and totally geodesic submanifold of $\HNC$, then $\Gamma$ is virtually Abelian.
			\item If $\Gamma$ is virtually Abelian, then $\HNC/\Gamma$ is a Stein manifold.
		\end{enumerate}
	\end{theoremIntro}
	
	The second point of Theorem \ref{thm-parabolique-abelien-delta-stein} is a consequence of Theorem \ref{thm-caracterisation-parab-stein}, stated in Section \ref{section-paraboliques}, which provides a complete characterisation of discrete and torsion-free parabolic subgroups $\Gamma$ of $\PU(n,1)$ for which $\HNC/\Gamma$ is Stein, and whose proof consists in reducing the problem to the unipotent case solved earlier in \cite{miebachQuotientsBoundedHomogeneous2024}. In complex dimension 2, this characterisation takes the following simpler form: a parabolic quotient $\mathbb{H}^{2}_{\C}/\Gamma$ is Stein if and only if $\Gamma$ is virtually Abelian (Corollary \ref{cor-parab-dim2}).

	In \cite{deyNoteComplexhyperbolicKleinian2020}, Dey and Kapovich have shown on the one hand that $\HNC/\Gamma$ does not contain any compact analytic subset of positive dimension if $\delta(\Gamma)<2$ and on the other hand that $\HNC/\Gamma$ is holomorphically convex if $\Gamma$ is convex-cocompact. By contrast, when $\Gamma$ is geometrically finite, the manifold $\HNC/\Gamma$ is not always holomorphically convex. For instance, when $\Gamma$ is a non-uniform lattice in $\PU(n,1)$, any holomorphic function $f:\HNC/\Gamma\to \C$ is constant. Using techniques from \cite{chenDiscreteGroupsHolomorphic2013}, we get the following characterization.
	
	\begin{restatable}{theoremIntro}{thmquotientholcvx}\label{thm-quotient-hol-cvx}
		Let $\Gamma$ be a geometrically finite and torsion-free subgroup of $\PU(n,1)$. The following are equivalent: \begin{enumerate}
			\item The manifold $\HNC/\Gamma$ admits a plurisubharmonic exhaustion function.
			\item For any parabolic subgroup $P<\Gamma$, the quotient $\HNC/P$ is holomorphically convex, or equivalently a Stein manifold.
			\item The manifold $\HNC/\Gamma$ is holomorphically convex.
		\end{enumerate}
	\end{restatable}
	
	 The strategy for proving Theorem \ref{thm-geom-fini-delta-stein} is to use Theorems \ref{thm-parabolique-abelien-delta-stein} and \ref{thm-quotient-hol-cvx} by showing that if $\Gamma$ satisfies one of the assumptions of Theorem \ref{thm-geom-fini-delta-stein}, then its parabolic subgroups are virtually Abelian. In connection with the conjecture of Dey and Kapovich, we also show that quotients $\HNC/\Gamma$ with $\delta(\Gamma)<2$ always admit non-constant holomorphic functions, as a particular case of Proposition \ref{prop-KH-ss-ens} below.
	
	We provide in Section \ref{section-paraboliques} an example showing that the constant 2 in Theorem \ref{thm-parabolique-abelien-delta-stein} is optimal, and an example of a unipotent parabolic group $\Gamma$ for which $\HNC/\Gamma$ is biholomorphic to a bundle of punctured disks over a non-compact Cousin manifold, and is not holomorphically convex. By a Cousin manifold, we mean a quotient of $\C^{n-1}$ by a discrete subgroup, which does not admit any non-constant holomorphic function, see \cite{cousinFonctionsTriplementPeriodiques1910,kopfermannMaximaleUntergruppenAbelscher1964,abeToroidalGroups2001}. The critical exponent of this example is equal to $\frac{5}{2}$.
	
	A \textit{complex Fuchsian group} is a discrete and torsion-free subgroup $\Gamma$ of $\PU(n,1)$ which acts cocompactly on a $\Gamma$-invariant complex geodesic. If $\Gamma$ is a complex Fuchsian group, then its critical exponent is equal to 2, and the quotient $\HNC/\Gamma$ contains a compact subvariety\footnote{Here and throughout all this article we use the word \textit{subvariety} as a synonym of \textit{closed analytic subset}.} of positive dimension. When $\Gamma$ is a convex-cocompact and torsion-free subgroup of $\PU(n,1)$ with critical exponent $\delta(\Gamma)=2$, Dey and Kapovich conjecture that $\HNC/\Gamma$ is non-Stein if and only if $\Gamma$ is a complex Fuchsian group, see \cite[Conjecture 17]{deyNoteComplexhyperbolicKleinian2020}. Using Patterson--Sullivan theory and techniques from \cite{connellNaturalFlowCritical2023}, we confirm this conjecture as follows.
	
	\begin{theoremIntro}\label{thm-delta2}
		Let $\Gamma$ be a discrete and torsion-free subgroup of $\PU(n,1)$ with $\delta(\Gamma)=2$. Suppose that $\HNC/\Gamma$ contains a compact subvariety of positive dimension. Then $\Gamma$ is a complex Fuchsian group.
	\end{theoremIntro}	
	
	We now discuss the relation of our results with earlier works. Corollary \ref{cor-free} applies in particular to representations of free groups seen as finite-index subgroups of the examples in \cite{goldmanComplexHyperbolicIdeal1992,falbelModuliSpaceModular2003}. It also applies to Schottky quotients, recovering \cite[Theorem 4.3]{miebachSchottkyGroupActions2018}. Theorem \ref{thm-parabolique-abelien-delta-stein}, together with the fact that the quotient of $\HNC$ by a loxodromic cyclic group is Stein (see for example \cite{defabritiisFamilyComplexManifolds1998manuel}, \cite{chenDiscreteGroupsHolomorphic2013} or Section \ref{section-kahler-hadamard} below), implies that the quotient of the complex hyperbolic space by any discrete and torsion-free Abelian subgroup of $\PU(n,1)$ is a Stein manifold. In \cite{chenDiscreteGroupsHolomorphic2013}, Chen asks whether the quotient of $\HNC$ by a discrete and torsion-free subgroup of $\PO(n,1)$ is Stein. Theorems \ref{thm-geom-fini-delta-stein} and \ref{thm-parabolique-abelien-delta-stein} yield a positive answer to this question for geometrically finite  or parabolic subgroups.

	Here are some earlier results about the analytic properties of quotients of the complex hyperbolic space $\HNC$. It is known that the quotient of $\HNC$ by an infinite discrete cyclic group is a Stein manifold \cite{defabritiisFamilyComplexManifolds1998manuel,defabritiisQuotientsUnitBall2001manuel,miebachQuotientsBoundedHomogeneous2010}. The article \cite{chenDiscreteGroupsHolomorphic2013} gives criteria for a quotient of $\HNC$ by a discrete subgroup to be Stein, and in particular shows that a quotient of the complex hyperbolic space by a unipotent Abelian parabolic group is Stein. It also contains results in the more general setting of quotients of Kähler-Hadamard manifolds. The case of quotients by unipotent parabolic subgroups is completely solved in \cite{miebachQuotientsBoundedHomogeneous2024}. Finally, as mentioned above, the article \cite{deyNoteComplexhyperbolicKleinian2020} in which the above conjecture appears contains the analogue of Theorem \ref{thm-geom-fini-delta-stein} in the case of convex-cocompact groups. Section 9 of the overview article \cite{kapovichSurveyComplexHyperbolic2022} also contains interesting results on the ends of $\HNC/\Gamma$, which are related to the analytic properties of this manifold.
	
	In another direction, we show that \cite[Theorem 1]{deyNoteComplexhyperbolicKleinian2020} can be extended to Kähler-Hadamard manifolds $(X,\omega)$ of negative curvature, using Patterson--Sullivan theory. We obtain the following result, where the critical exponent of $\Gamma$ is defined by the same Formula \eqref{eq-def-delta} as for discrete subgroups of $\PU(n,1)$, by choosing an arbitrary point $o$ of $X$ and with $d$ the Riemannian distance associated with $\omega$. This number does not depend on the choice of $o\in X$.
	
	\begin{propositionIntro}\label{prop-KH-general}
		Let $(X,\omega)$ be a simply connected complete Kähler manifold with sectional curvature bounded above by $-1$, and $\Gamma$ be a group acting freely and properly discontinuously by holomorphic isometries on $X$.\begin{enumerate}
			\item If $\Gamma$ acts cocompactly on a $\Gamma$-invariant convex subset of $X$, then $X/\Gamma$ is holomorphically convex.
			\item If the critical exponent of $\Gamma$ is less than two, then $X/\Gamma$ does not contain any compact subvariety of positive dimension. 
		\end{enumerate}
	\end{propositionIntro}
	
	\noindent More precise results that imply Proposition \ref{prop-KH-general} are given in Subsection \ref{ss-section-preuve-prop5}.
	
	This article is organised as follows. In the first section are recalled the notions of holomorphic convexity, Stein manifolds, geometrically finite groups and Patterson--Sullivan measures, followed by the proof of Proposition \ref{prop-KH-general}. In Subsection \ref{ss-section-criteres}, we give criteria for asserting that a quotient of the form $X/\Gamma$ does not admit any compact subvariety of positive dimension, with $X$ and $\Gamma$ as above. Theorem \ref{thm-parabolique-abelien-delta-stein} is proved in Section \ref{section-paraboliques}. Section \ref{section-hol-cvx-geom-fini} contains the proof of Theorem \ref{thm-quotient-hol-cvx} from which is deduced the proof of Theorem \ref{thm-geom-fini-delta-stein} and Corollary \ref{cor-free}. Section \ref{section-delta2} is independant from Sections \ref{section-paraboliques} and \ref{section-hol-cvx-geom-fini}, and contains a proof of Theorem \ref{thm-delta2}.
	
	\paragraph{Acknowledgments.} I thank Pierre Py for his constant support and help throughout this work and Christian Miebach for interesting conversations we had in Calais. I also thank the referees for their comments and suggestions on the text.
	
	\tableofcontents
	
	\section{Quotients of negatively curved Kähler-Hadamard manifolds}\label{section-kahler-hadamard}

	We begin by recalling some definitions from complex analysis and then review some results of negative curvature geometry. After that, we state and prove results that imply Proposition \ref{prop-KH-general} in Subsection \ref{ss-section-preuve-prop5}. Finally, we summarize in Subsection \ref{ss-section-criteres} some known criteria for asserting that $X/\Gamma$ does not admit a compact subvariety of positive dimension.
	
	\subsection{Generalities on Stein manifolds}\label{ss-section-stein}
	
	We first recall the definition of plurisubharmonic and strictly plurisubharmonic functions, and we refer to \cite[Chapter 1]{demaillyComplexAnalyticDifferential-manuel} for more details. Let $X$ be a complex manifold. A continuous function $f:X\to \R$ is \textit{plurisubharmonic} if for every chart $\phi:V\subset X\longrightarrow W\subset \C^{n}$, every $a\in W$ and every $\xi\in \C^{n}$ such that $\lVert \xi\rVert < d(a,{}^{c}W)$, we have \begin{equation*}
		f\circ\phi^{-1}(a) \le \frac{1}{2\pi}\int_{0}^{2\pi}f\circ\phi^{-1}(a+e^{i\theta}\xi)d\theta.
	\end{equation*}
	It is \textit{strictly plurisubharmonic} if for every $x\in X$ there are holomorphic coordinates $(z_1,\dots,z_n)$ defined on some neighborhood of $x$ and a constant $c>0$ such that $z\mapsto f(z)-c\lVert z\rVert^{2}$ is plurisubharmonic. If $f$ is of class $\mathcal{C}^2$, then it is plurisubharmonic (resp. strictly plurisubharmonic) if and only if the $(1,1)$-form $i\partial\bar\partial f$ is nonnegative (resp. positive).
	
	A complex manifold $X$ is said to be \textit{holomorphically convex} if the \textit{holomorphic hull} $\widehat{K}$ of any compact subset $K$ of $X$, which is defined by\begin{equation*}
		\widehat{K} := \{x\in X \ \vert\ \forall f\in \mathcal{O}(X), \lvert f(x)\rvert \le \sup_{K}\lvert f\rvert\},
	\end{equation*}
	is compact. The manifold $X$ is \textit{Stein} if it is holomorphically convex and if, in addition, for any pair $(x,y)$ of distinct points of $X$, there is a holomorphic function $f:X\longrightarrow \C$ such that $f(x)\ne f(y)$. Grauert's theorem asserts that a manifold is Stein if and only if it admits a strictly plurisubharmonic exhaustion function, see \cite{grauertLeviProblemImbedding1958}. Alternatively, a manifold is Stein if and only if it is holomorphically convex and does not contain any compact subvariety of positive dimension. This follows from the existence of the Remmert reduction of a holomorphically convex manifold \cite[Theorem 2.1]{peternellPseudoconvexityLeviProblem1994}. We will also use the following result, that we subsequently refer to as Grauert's theorem, since it derives from it.
	
	\begin{theorem*}[{\cite{grauertLeviProblemImbedding1958}, \cite[Corollary 2.4]{peternellPseudoconvexityLeviProblem1994}},\cite{narasimhanLeviProblemComplex1962}]
		Let $X$ be a complex manifold admitting a continuous exhaustion function which is strictly plurisubharmonic outside a compact set. Then $X$ is holomorphically convex.
	\end{theorem*} 
	
	We will also use the following classical result.
	
	\begin{theorem*}\label{thm-docquier-grauert}
		Let $X$ be a complex manifold and $f:X\to \R$ be a strictly plurisubharmonic continuous function. Assume that for every real number $t$, the open subset \begin{equation*}
			X_t := \{x\in X\mid f(x)<t\}
		\end{equation*}
		of $X$ is a Stein manifold. Then $X$ is a Stein manifold.
	\end{theorem*}
	
	\begin{proof}
		For every integer $n$, the open set $X_{n+1}$ is Stein and by \cite[Corollary 1]{narasimhanLeviProblemComplex1962}, we get that $(X_{n},X_{n+1})$ is a Runge pair. Using \cite[Satz 1.3]{steinUberlagerungenHolomorphvollstandigerKomplexer1956}, we deduce that $X$ is Stein.
	\end{proof}

	\subsection{Convexity and Busemann functions on Kähler-Hadamard manifolds}

	In this subsection, $(X,\omega)$ denotes a simply connected complete Kähler manifold with complex structure denoted by $J$ and sectional curvature bounded above by $-1$.
	
	Let $d$ be the Riemannian distance associated with $\omega$. We recall that if $\phi:X\to \R$ is a function of class $\mathcal{C}^2$, then the form $i\partial\overline{\partial}\phi$ is related to the Riemannian Hessian $D^{2}(\phi)$ of $\phi$ by  \begin{equation}\label{eq-lev-hessienne}
		2i\partial\overline{\partial}\phi(v,Jv)=D^{2}(\phi)(v,v)+D^{2}(\phi)(Jv,Jv)
	\end{equation}
	for all tangent vector $v$ of $X$, see \cite{greeneSubharmonicityPlurisubharmonicityGeodesically1973}.
	
	 A continuous function $\phi:X\longrightarrow \R$ is called \textit{convex} if for all geodesics $\eta:\R\to X$, the function $\phi\circ\eta$ is convex, and it is called \textit{strictly convex} if for any compact subset $K\subset X$, there exists a constant $\alpha>0$ such that, for any unit-speed geodesic $\eta:[0,1]\to K$, the function $t\in [0,1]\mapsto \phi\circ\eta(t)-\alpha t^{2}$ is convex. Since $X$ is a Kähler manifold, every continuous convex function $X\to R$ is plurisubharmonic \cite[Theorem 3]{greeneSubharmonicityPlurisubharmonicityGeodesically1973}. It follows that a strictly convex function $f:X\to \R$ is strictly plurisubharmonic. Indeed, let $x$ be a point in $X$ and $(z_1,\dots,z_n)$ be holomorphic coordinates defined in a neighborhood $\Omega$ of $x$. Then for any open subset $V$ with $\overline{V}\subset \Omega$, there exists a constant $c>0$ such that $z\mapsto f(z)-c\lVert z\rVert^{2}$ is convex, and therefore plurisubharmonic, on $V$. Thus $f$ is strictly plurisubharmonic. When $f$ is of class $\mathcal{C}^2$, these statements are a simple consequence of Formula \eqref{eq-lev-hessienne}.
	
	We denote by $\partial X$ the visual boundary of $X$, defined as the quotient of the set of geodesic rays in $X$ by the equivalence relation ``remaining at bounded distance''. For every point $p\in X$, the unit tangent space at $p$ is in bijection with $\partial X$, and we endow $\partial X$ with the topology that makes this bijection a homeomorphism. Any isometry of $X$ induces a homeomorphism of $\partial X$.
	
	Fix a point $o\in X$. For all $\xi \in \partial X$, the \textit{Busemann function} at $\xi$ is the function defined by\begin{equation*}
		\forall x\in X, B_{\xi}(x) := B_{\xi}(x,o) := \lim_{z\to \xi}(d(x,z)-d(o,z)).
	\end{equation*}
	The Busemann function at $\xi$ depends on $o\in X$ only up to an additive constant. This function is of class $\mathcal{C}^{2}$, see \cite[Proposition 3.1]{heintzeGeometryHorospheres1977}, and it depends continuously on $\xi \in \partial X$. It is moreover strictly plurisubharmonic and more precisely we have\begin{equation}\label{eq-Busemann-spsh}
		i\partial\overline{\partial}B_{\xi} \ge \omega.
	\end{equation}
	This inequality is a consequence of \cite[Proposition 2.28]{greeneFunctionTheoryManifolds1979}, see also \cite{siuCompactificationNegativelyCurved1982,chenDiscreteGroupsHolomorphic2013}. We also call $(\xi,x,y)\in \partial X\times X^{2} \longmapsto B_{\xi}(x,y) \in\R$ the Busemann function on $X$.
	
	The gradient of $B_{\xi}$ at a point $x\in \HNC$ is $-v_{x\xi}$, where $v_{x\xi}$ is the unit tangent vector at $x$ pointing at $\xi$. In particular, if we denote by $h$ the Hermitian product associated with $\omega$, then for all $\xi\in \partial X$ and for all $v\in T_xX$, we have\begin{equation}\label{ineq-dBJ}
		dB_{\xi}(v)^{2} + dB_{\xi}(Jv)^{2} = \lvert h(v,v_{x,\xi})\rvert^{2}\le \lVert v\rVert^{2}.
	\end{equation}
A sublevel set of $B_{\xi}$ is called a \textit{horoball} at $\xi$.	
\subsection{Discrete subgroups, Patterson--Sullivan measure and geometrical finiteness}\label{sec:gen-geom-finite}
	
	As in the previous subsection $(X,\omega)$ is a simply connected complete Kähler manifold with complex structure denoted by $J$ and sectional curvature bounded above by $-1$. In addition, let $\Gamma$ be a group acting freely and properly discontinuously by holomorphic isometries on $X$. We also assume that $\Gamma$ is non-elementary, which means that $\Gamma$ does not stabilize a geodesic of $X$, nor a point of $\partial X$.\par

The \textit{limit set} $\Lambda(\Gamma)$ of $\Gamma$ is the closed subset of $\partial X$ defined as the accumulation set of an orbit $\Gamma \cdot o$, for some point $o\in X$, and it does not depend on the choice of the point. The \textit{domain of discontinuity} $\Omega(\Gamma)$ of $\Gamma$ is an open subset of $\partial X$ which can be defined as the complement of the limit set. These sets are invariant under the action of $\Gamma$ on $\partial X$, and in particular, the closed geodesic convex hull of the limit set forms a $\Gamma$-invariant closed subset of $X$. The quotient $C_{\Gamma}$ of this convex hull by $\Gamma$ is a subset of $X/\Gamma$, called \textit{convex core} of $X/\Gamma$.

To define the notion of geometrical finiteness, we need the notions of conical limit points and bounded parabolic points. First, a point $\xi\in \partial X$ is called \textit{conical} if it is the limit of a sequence in $(\Gamma \cdot o)^{\N}$ that stays at bounded distance of a(ny) geodesic ray pointing to $\xi$. The point $\xi$ is \textit{parabolic} if it is the fixed point of a parabolic element of $\Gamma$. It is a \textit{bounded parabolic point} if it is parabolic and if the action of the parabolic subgroup $\operatorname{Stab}_{\Gamma}(\xi)$ on $\Lambda(\Gamma)\setminus\{\xi\}$ is cocompact.

\begin{definition}[{\cite[Definition F2]{bowditchGeometricalFinitenessVariable1995}}]
	The group $\Gamma$ is called \textit{geometrically finite} if every element $\xi\in\Lambda(\Gamma)$ is either a bounded parabolic point or a conical point.
\end{definition}

The structure of geometrically finite quotients of $X$ is described by the following theorem, which relies on the \textit{thick-thin decomposition} and on \textit{Margulis lemma}. We refer for instance to \cite[Chapter D]{benedettiLecturesHyperbolicGeometry1992} or \cite[Section 3.5]{bowditchGeometricalFinitenessVariable1995} for an account of these two notions.

\begin{theorem*}[{\cite[Definition F1]{bowditchGeometricalFinitenessVariable1995}}]
	The group $\Gamma$ is geometrically finite if and only if for a(ny) positive $\epsilon$ less than the Margulis constant $\epsilon_0$, the intersection of the convex core and the $\epsilon$-thick part of $X/\Gamma$ is compact.
	
	Moreover, in this case, the $\epsilon$-thick part $Q$ of $X/\Gamma$ is relatively compact in the quotient $\left(X\cup \Omega(\Gamma)\right)/\Gamma$, and the thin part of $X/\Gamma$ consists of a finite number of parabolic ends, meaning that $X_{\Gamma}:=X/\Gamma$ decomposes as 
	\begin{equation}\label{eq-decompo-GF}
		X_{\Gamma} =: Q \cup \bigcup_{i=1}^{k}E_i,
	\end{equation}
	where $k$ is an integer, and for $i\in \{1,\dots,k\}$, each $E_i$ is an open subset of $X_{\Gamma}$ biholomorphic to the quotient of a horoball $B_i^{-1}((-\infty,0))$ by a maximal parabolic subgroup $P_i$ of $\Gamma$, for some Busemann function $B_i$.
\end{theorem*}
When $X$ is a hyperbolic space (either real, complex, quaternionic or octonionic) and $\Gamma$ is a geometrically finite group acting by isometries on $X$, we recall that the critical exponent of $\Gamma$ equals the Hausdorff dimension of its limit set $\Lambda(\Gamma)$ for an appropriate distance defined on $\partial X$, see \cite[Theorem 6.1]{corletteLimitSetsDiscrete1999}. This result will not be used in the sequel of the article.

In this article, we will also use Patterson-Sullivan theory. It is used in the proofs of Proposition \ref{prop-KH-ss-ens}, Corollary \ref{cor:delta-croissance-ssexp}, and Theorem \ref{thm-delta2}. We now recall the definition and some basic properties of Patterson-Sullivan measures, and we refer the reader to \cite{roblinErgodiciteEquidistributionCourbure2003,pattersonLimitSetFuchsian1976, sullivanDensityInfinityDiscrete1979, nichollsErgodicTheoryDiscrete1989} for more details and for the construction of Patterson--Sullivan measures. Let $X$ and $\Gamma$ be as above, and $\delta$ be the critical exponent of $\Gamma$, whose definition was recalled in the introduction. A Patterson-Sullivan measure is a $\Gamma$-conformal density of dimension $\delta$, which means that it is a family of finite measures $(\mu_x)_{x\in X}$ on $\partial X$ such that $\gamma_*\mu_x = \mu_{\gamma x}$ for all $x\in X$ and $\gamma\in \Gamma$, and such that
		\begin{equation}\label{eq-prop-PS} \forall x,y\in X, \frac{d\mu_x}{d\mu_y}=e^{-\delta B_{\bullet}(x,y)}. \end{equation} 
		Moreover, the support of the measure $\mu_x$ is the limit set $\Lambda(\Gamma)$ of $\Gamma$.

\subsection{Proof of Proposition \ref{prop-KH-general}}\label{ss-section-preuve-prop5}

We now state results that are more precise than Proposition \ref{prop-KH-general} and that imply it. In the following two propositions and the subsequent corollaries, which provide an alternative proof of the results of \cite{deyNoteComplexhyperbolicKleinian2020}, $(X,\omega)$ and $\Gamma$ are as in the previous subsection. In this context, we mean by \textit{pinched}, when referring to the sectional curvature of $(X,\omega)$, that it is bounded below by $-b$, and when referring to the Ricci curvature, that it is bounded below by $-b\omega(\cdot,J\cdot)$, for some constant $b>1$.

\begin{proposition}\label{prop-KH-hol-cvx}
	Let $C$ be a $\Gamma$-invariant closed and geodesically convex subset of $X$. Then the compact connected subvarieties of positive dimension of $X/\Gamma$ are included in $C/\Gamma$. Moreover, if the action of $\Gamma$ on $C$ is cocompact, then $X/\Gamma$ is holomorphically convex.
\end{proposition}

\begin{proposition}\label{prop-KH-ss-ens}
	Assume that $\delta(\Gamma)<2$. Then $X/\Gamma$ admits a strictly plurisubharmonic function and in particular does not contain any compact subvariety of positive dimension. If moreover $X$ has pinched Ricci curvature, then holomorphic functions on $X/\Gamma$ separate points and define local coordinates at all points of $X/\Gamma$.
\end{proposition}

\begin{corollary}\label{cor-cvx+delta:stein}
	Assume that $\delta(\Gamma)<2$ and $\Gamma$ is convex-cocompact. Then $X/\Gamma$ is a Stein manifold.
\end{corollary}

Our next corollary involves an assumption on the Patterson--Sullivan measure associated with $\Gamma$. The Patterson--Sullivan measure $(\mu_x)_{x\in X}$ is said to have subexponential growth if for all $\eta>0$, there is a constant $C_{\eta}>0$ such that\begin{equation*}
	\forall x\in X,\  \mu_x(\partial X) \le C_{\eta}e^{\eta d(x,o)},
\end{equation*}
for some basepoint $o\in X$, see \cite[§1.4]{connellNaturalFlowCritical2023}. For instance, if $X$ has pinched sectional curvature, $X/\Gamma$ has positive injectivity radius and if the Bowen-Margulis measure $\nu$ associated with $(\mu_x)_{x\in X}$ is finite, then the total masses $\mu_x(\partial X)$ are uniformly bounded, see \cite[Theorem 1.15]{connellNaturalFlowCritical2023}. We also refer to the latter article for a definition of the Bowen-Margulis measure $\nu$.

\begin{corollary}\label{cor:delta-croissance-ssexp}
	Assume that $\delta(\Gamma)<2$, $X$ has pinched sectional curvature and the Patterson-Sullivan measure $(\mu_x)_{x\in X}$ has subexponential growth. Then $X/\Gamma$ is a Stein manifold.
\end{corollary}

Proposition \ref{prop-KH-hol-cvx} is proven by an application of the next lemma to the square of the distance function to $C/\Gamma\subset\HNC/\Gamma$.

\begin{lemma}\label{lem-cvx-psh}
	Let $(M,\omega)$ be a Kähler manifold. Assume that there exists a continuous function $\phi:M\to \R$ which is convex on $M$, and strictly convex on $M\setminus \phi^{-1}(0)$. Then the compact connected subvarieties of positive dimension of $M$ are included in the level set $\phi^{-1}(0)$. If moreover $\phi$ is an exhaustion, then $M$ is holomorphically convex.
\end{lemma}

\begin{proof}
	Let $\phi:M\to \R$ be a continuous function which is convex on $M$, and strictly convex on $M\setminus \phi^{-1}(0)$. Then $\phi$ is plurisubharmonic on $M$ and strictly plurisubharmonic on $M\setminus \phi^{-1}(0)$.  Moreover for any connected compact subvariety $A$ of $M$, the function $\phi\vert_{A}$ is constant by the maximum principle, so $A$ is included in $\phi^{-1}(0)$ or in $M\setminus \phi^{-1}(0)$. Notice that $A$ cannot be contained in $M\setminus \phi^{-1}(0)$ because in that case $\phi\vert_{A}$ would be constant and strictly plurisubharmonic. Thus $A\subset \phi^{-1}(0)$. If moreover $\phi$ is an exhaustion, then $M$ is holomorphically convex by Grauert's theorem.
\end{proof}

	\begin{proof}[Proof of Proposition \ref{prop-KH-hol-cvx}]
		Let $d_C^{2}:X\to \R_+$ be the square of the distance function to $C$. This function is convex on $X$, and it is strictly convex on $X\setminus C$. This is proved in \cite[Lemma 4.5]{benoistHarmonicQuasiisometricMaps2023}, and for the sake of completeness we now outline a proof of the strict convexity of $d_C^{2}$ on $X\setminus C$. Fix some $\epsilon>0$ and let $\gamma:[0,L]\to X$ be a unit-speed geodesic such that $d(\gamma(0),C)\ge \epsilon$ and $d(\gamma(L),C)\ge \epsilon$. Denote by $x$ and $y$ the projections of $\gamma(0)$ and $\gamma(L)$ in $C$ and let $\eta:[0,L']\to X$ be the unit-speed geodesic joining $x$ and $y$. Using \cite[Corollary 2.1.3 - Formula 2.1(iv)]{korevaarSobolevSpacesHarmonic1993}, we get :\begin{equation*}
			2d_C^{2}\left(\gamma\left(\frac{L}{2}\right)\right) - d_C^{2}(\gamma(0)) - d_C^{2}(\gamma(L)) \le -\frac{1}{2}(L-L')^{2}.
		\end{equation*}
		Finally, since $X$ has sectional curvature bounded above by $-1$, there is a positive constant $a$ depending only on $\epsilon$ such that $(L-L')^{2}\ge aL^{2}$. Using \cite[Lemma 1]{greeneConvexFunctionsManifolds1976manuel}, we conclude that $d_C^{2}$ is strictly convex on $X\setminus C$.
		
		To conclude, notice that the function $d_C^{2}$ is $\Gamma$-invariant, so it defines a convex function $\phi:X/\Gamma\to \R$, which is strictly convex outside $C/\Gamma$. The proposition is obtained by applying Lemma \ref{lem-cvx-psh} with the function $\phi$.
	\end{proof}
	
	\begin{remark}
		Proposition \ref{prop-KH-hol-cvx} provides an alternative proof of \cite[Proposition 5]{kapovichSurveyComplexHyperbolic2022}: suppose that $X$ has negatively pinched curvature and that there is a surjective holomorphic map $f:X/\Gamma \to B$ with compact fibers on a complex manifold $B$ with $\dim(B)<\dim(X)$. Then $\Lambda(\Gamma)=\partial X$. In particular $X/\Gamma$ cannot have convex ends. 
	\end{remark}

	The next lemma, used in the proof of Proposition \ref{prop-KH-ss-ens}, asserts that a certain function defined in \cite{connellNaturalFlowCritical2023} is strictly plurisubharmonic when the critical exponent of $\Gamma$ is less than 2. This function and the flow it defines also play an important role in Section \ref{section-delta2}. An alternative proof of the first point of Proposition \ref{prop-KH-ss-ens}, which uses comparison arguments from \cite{greeneFunctionTheoryManifolds1979} is outlined below. 
	
	\begin{lemma}\label{lem-f-unif-spsh}
		Let $(X,\omega)$ be a simply connected complete Kähler manifold with sectional curvature bounded above by $-1$, and $\Gamma$ be a non-elementary group acting freely and properly discontinuously by holomorphic isometries on $X$. Denote by $\delta$ the critical exponent of $\Gamma$, by $(\mu_x)_{x\in X}$ a Patterson-Sullivan measure associated with $\Gamma$ and by $\lVert \mu_x\rVert$ the total mass of the measure $\mu_x$ for every $x\in X$. Then the $\Gamma$-invariant function on $X$ defined by $f(x):=-\ln \lVert \mu_x\rVert$ satisfies \begin{equation*}
			i\partial\overline{\partial}{f} \ge \delta(1-\frac{\delta}{2})\omega.
		\end{equation*}
	\end{lemma}
	
	\begin{proof}
		For every $x\in X$, let $\overline{\mu_x}$ be the normalized probability measure $\overline{\mu_x} = \frac{\mu_x}{\lVert \mu_x\rVert}$. Fixing a point $o\in X$, denote by $B_\theta:=B_{\theta}(\cdot,o)$ the Busemann function at a point  $\theta\in \partial X$ which vanishes at $o$. From dominated convergence together with the $\mathcal{C}^2$-regularity of Busemann functions and Formula \eqref{eq-prop-PS}, it follows that $f$ is of class $\mathcal{C}^2$.
		We will now compute the Levi form of $f$. First, the differential of $f$ at a point $x\in X$ is given by
		\begin{equation*}
			df(x) = \frac{\int_{\partial X}\delta dB_{\theta}(x)e^{-\delta B_{\theta}(x)}d\mu_o(\theta)}{\lVert \mu_x\rVert}.
		\end{equation*}
		Then, its Hessian is computed as follows.
		\begin{align*}
			\begin{split}
				D^{2}f(x) &= \frac{1}{\lVert \mu_x\rVert}\int_{\partial X}\left(\delta D^{2}B_{\theta}(x)e^{-\delta B_{\theta}(x)}-\delta^{2}dB_{\theta}(x)\otimes dB_{\theta}(x)e^{-\delta B_{\theta}(x)}\right)d\mu_o(\theta) +\\ & \hspace{5.6cm} \frac{\delta^{2}}{\lVert \mu_x\rVert^{2}}\left(\int_{\partial X}dB_{\theta}(x)e^{-\delta B_{\theta}(x)}d\mu_o(\theta)\right)^{2}\\
				&={\delta}\int_{\partial X}D^{2}B_{\theta}(x)d\overline{\mu_x}(\theta) +\\& \hspace{2.4cm} {\delta^{2}}\left(\left(\int_{\partial X}dB_{\theta}(x)d\overline{\mu_x}(\theta)\right)^{2}-\int_{\partial X}dB_{\theta}(x)\otimes dB_{\theta}(x)d\overline{\mu_x}(\theta)\right).
			\end{split}
		\end{align*}
	Let $v\in T_xX$. Using Identity \eqref{eq-lev-hessienne}, we obtain
	\begin{equation*}
		i\partial{\partial}\overline{f}(v,Jv) \ge \delta \int_{\partial X}i\partial\overline{\partial}B_{\theta}(v,Jv)d\overline{\mu_x}(\theta) - \frac{\delta^{2}}{2}\int_{\partial X}(dB_{\theta}(v)^{2}+dB_{\theta}(Jv)^{2}) d\overline{\mu_x}(\theta).
	\end{equation*}	
		Using Inequalities \eqref{eq-Busemann-spsh} and \eqref{ineq-dBJ}, we deduce that
		\begin{equation*}
			i\partial\overline{\partial}{f} \ge \delta(1-\frac{\delta}{2})\omega.\qedhere
		\end{equation*}
	\end{proof}
	
	\begin{proof}[Proof of Proposition \ref{prop-KH-ss-ens}]
		If $\Gamma$ is elementary, the result is already known, see \cite[Theorem 1.1 and Proposition 1.3]{chenDiscreteGroupsHolomorphic2013}. Assume that $\Gamma$ is non-elementary and that $\delta<2$. Then the function $f$ defined in Lemma \ref{lem-f-unif-spsh} is strictly plurisubharmonic. Suppose moreover that there is a constant $C>0$ such that\begin{equation*}
			\Ricci(\omega)\ge -C\omega(\cdot,J\cdot).
		\end{equation*}
		Then there is a constant $C'>0$ such that\begin{equation*}
			i\partial\overline{\partial}(C'f) + \Ricci(\omega)(\cdot,J\cdot) \ge 0.
		\end{equation*}
		Using \cite[Proposition 4.1]{chenDiscreteGroupsHolomorphic2013}, we obtain that the holomorphic functions on $X/\Gamma$ separate points and give local coordinate systems.
	\end{proof}
	
	\begin{remark}A strictly plurisubharmonic function on $X/\Gamma$ can also be constructed following Dey and Kapovich's ideas, providing an alternative proof for the first point of Proposition \ref{prop-KH-ss-ens}. Here is an outline of the argument. Let $\phi$ be the function defined on $X$ by $\phi(x) :=\tanh(d(o,x))^{2}$ for some basepoint $o\in X$. An application of the comparison result \cite[Theorem A]{greeneFunctionTheoryManifolds1979} together with Formula \eqref{eq-lev-hessienne} gives that $\phi$ is strictly plurisubharmonic on $X$: this is obtained by comparing the Hessian of $\phi$ with the Hessian of the function $\widetilde{\phi}$ defined on $\mathbb{H}^{2n}_{\R}$, the real hyperbolic space of dimension $2n$, by $\widetilde{\phi}(x) := \tanh(d_{\text{ hyp}}(\widetilde{o},x))^{2}$ for some basepoint $\widetilde{o}$ of $\mathbb{H}^{2n}_{\R}$. Then because\begin{align*}
				0\leq 1-\phi \leq 4e^{-2d(o,\cdot)},
			\end{align*}
			we deduce that the convergence of the series \begin{equation*}
				\sum_{\gamma \in \Gamma}e^{-2d(o,\gamma o)}
			\end{equation*} implies that \begin{equation*}
				\sum_{\gamma \in \Gamma}(\phi(\gamma\cdot x)-1)
			\end{equation*} converges uniformly on compact subsets of $X$ to a strictly plurisubharmonic and $\Gamma$-invariant function $\psi$. When $X=\HNC$, the function $\phi$ is the squared euclidean norm on the unit ball, and the above series is the one constructed in \cite{deyNoteComplexhyperbolicKleinian2020}.
		\end{remark}

\begin{proof}[Proof of Corollary \ref{cor-cvx+delta:stein}]
	This follows directly from Propositions \ref{prop-KH-hol-cvx} and \ref{prop-KH-ss-ens}, using that a manifold is Stein if and only if it is holomorphically convex and does not contain any compact subvariety of positive dimension. Alternatively, each one of the three strictly plurisubharmonic functions $\phi+e^{f}$, $\phi+e^{\psi}$ and $f$ is an exhaustion on $X/\Gamma$, with $\phi$ as in the proof of Proposition \ref{prop-KH-hol-cvx}, $f$ as in the proof of Proposition \ref{prop-KH-ss-ens} and $\psi$ as in the remark above. For the function $f$, this can be shown using that $(X\cup\Omega(\Gamma))/\Gamma$ is compact, where $\Omega(\Gamma):= \partial X\setminus \Lambda(\Gamma)$ is the discontinuity subset of $\Gamma$.
\end{proof}

\begin{proof}[Proof of Corollary \ref{cor:delta-croissance-ssexp}]
	Let $f$ be the $\Gamma$-invariant function defined in Lemma \ref{lem-f-unif-spsh} and $o$ be a point of $X/\Gamma$. By \cite{tamExhaustionFunctionsComplete2010} or \cite[§4 of Ch. 26]{chowRicciFlowTechniques2010}, there is a constant $C>0$ and a smooth function $g:X/\Gamma \to \R$ such that $g\ge d(\cdot,o)$ and $D^{2}(g)\ge -C\omega(\cdot,J\cdot)$. Using Identity \eqref{eq-lev-hessienne}, we deduce that $i\partial\overline{\partial}g\ge -C\omega$. Now define\begin{equation*}
		\phi = g + C'f,
	\end{equation*}
	where $C'$ is a positive constant such that\begin{equation*}
		-C+C'\delta(1-\frac{\delta}{2})>0.
	\end{equation*}
	Then $\phi$ is strictly plurisubharmonic by Lemma \ref{lem-f-unif-spsh}. Since $(\mu_x)_{x\in X}$ has subexponential growth, there is a constant $C''>0$ such that\begin{equation*}
		f\ge -\frac{1}{2C'}d(\cdot,o) -\ln C''.
	\end{equation*}
	We deduce that \begin{equation*}
		\phi \ge \frac{1}{2}d(\cdot,o) -\ln C'',
	\end{equation*}
	which implies that $\phi$ is proper. Thus $X/\Gamma$ is a Stein manifold.
\end{proof}
	
	\subsection{Compact subvarieties of positive dimension of $X/\Gamma$}\label{ss-section-criteres}
	We now summarize known criteria for asserting that $X/\Gamma$ does not contain a compact subvariety of positive dimension. 
	
	\begin{proposition}\label{prop-criteres}
		Let $(X,\omega)$ be a simply connected complete Kähler manifold with sectional curvature bounded above by $-1$, and $\Gamma$ be a group acting freely and properly discontinuously by holomorphic isometries on $X$. Denote by $d$ the Riemannian distance associated with $\omega$. Then each one of the following condition is sufficient to assert that $X/\Gamma$ does not contain a compact subvariety of positive dimension.\begin{enumerate}[label=(\alph*)]
			\item The group $\Gamma$ is parabolic in the sense that all its elements are parabolic isometries fixing the same point in $\partial X$.
			\item The critical exponent of $\Gamma$ satisfies $\delta(\Gamma)<2$.
			\item There exists a $\Gamma$-invariant geodesically convex subset $C$ of $X$ which is included in a totally real submanifold $M$ of $X$.
			\item The Kähler form $\omega$ is exact on $X/\Gamma$. 
			\item The cohomology group $H^{2}(\Gamma,\R)$ vanishes.
			\item There is a complete vector field on $X/\Gamma$ whose flow contracts complex subspaces. 
		\end{enumerate}
	\end{proposition}
	
	\begin{proof}\begin{enumerate}[leftmargin=*,label=(\alph*)]
			\item If $\Gamma$ is parabolic, there is a Busemann function at some point $\xi\in \partial X$ which is invariant under the action of $\Gamma$, see \cite[Proposition 7.8]{eberleinVisibilityManifolds1973}. As a consequence $X/\Gamma$ admits a strictly plurisubharmonic function.
			\item This follows from \cite{deyNoteComplexhyperbolicKleinian2020} for the complex hyperbolic space and from Proposition \ref{prop-KH-ss-ens} in the general case.
			\item This is a consequence of Proposition \ref{prop-KH-hol-cvx}. Indeed, let $A$ be a compact connected subvariety of positive dimension included in $X/\Gamma$. Then $A\subset C/\Gamma$, and taking a smooth point $x$ of the lift $\tilde{A}$ of $A$ in $X$, we obtain that $T_x\tilde{A}\cap J(T_x\tilde{A}) \subset T_xM$, which contradicts the hypothesis that $M$ is totally real. Alternatively, it  can be shown that in this case the distance squared function to the convex core is strictly plurisubharmonic, see \cite{chenDiscreteGroupsHolomorphic2013}.
			\item This is a well known fact about Kähler geometry.
			\item The manifold $X/\Gamma$ is a $K(\Gamma,1)$, and in particular its cohomology identifies with that of $\Gamma$. Therefore, if $H^{2}(\Gamma,\R)=0$, the Kähler form $\omega$ is exact on $X/\Gamma$. This argument appears in a dual version in \cite[Lemma 4.2]{miebachSchottkyGroupActions2018}, and also in \cite[page 27]{kapovichSurveyComplexHyperbolic2022}.
			\item See for instance \cite[Proof of Theorem 1.5]{connellNaturalFlowCritical2023}.\qedhere
		\end{enumerate}
	\end{proof}
	We now give examples of quotients $X/\Gamma$ that do not contain a compact subvariety of positive dimension. In these examples $X=\HNC$, the complex hyperbolic space of dimension $n$, whose construction and elementary properties are recalled in Subsection \ref{ss-section-preli-parab} below.
	\begin{examples}\begin{enumerate}
			\item When $X=\HNC$ and $\Gamma$ is a discrete and torsion-free subgroup of $\PO(n,1)$, embedded in $\PU(n,1)$ so as to stabilize a copy of $\HNR$, the quotient $\HNC/\Gamma$ does not contain a compact subvariety of positive dimension by Proposition \ref{prop-criteres}, c). In particular, as proven in \cite{chenDiscreteGroupsHolomorphic2013}, the quotient of $X$ by the group generated by a hyperbolic element is a Stein manifold.
			\item Let $\Gamma=\pi_1(\Sigma_g)$ be a surface group and $\rho:\Gamma \to \PU(n,1)$ be a discrete and faithful representation of $\Gamma$ in $\PU(n,1)$. The \textit{Toledo invariant} $\tau$ of $\rho$ is the real number \begin{equation*}
				\tau := \frac{1}{2\pi}[\phi^{*}\omega] \in H^{2}(\Sigma_g,\R) \simeq \R,
			\end{equation*}
			where $\phi : \Sigma_g \to \HNC/\rho(\Gamma)$ is any homotopy equivalence between $\Sigma_g$ and $\HNC/\rho(\Gamma)$. Then $\omega$ is exact on $\HNC/\rho(\Gamma)$ if and only if $\tau=0$. Thus representations of surface groups with $\tau=0$ provide examples of discrete subgroups $\Gamma$ of $\PU(n,1)$ such that $\HNC/\Gamma$ does not contain a compact subvariety of positive dimension. These examples appear in \cite[page 27]{kapovichSurveyComplexHyperbolic2022}.
			\item If $\Gamma$ is a free group, then $H^{2}(\Gamma,\R)=0$ so $X/\Gamma$ does not contain any compact subvariety of positive dimension.
			\item Let $\Gamma$ be a uniform and torsion-free lattice in $\PU(n,1)$. There is a holomorphic map  $A$ from $\HNC/\Gamma$ to an Abelian variety of dimension $N=\frac{b_1(\Gamma)}{2}$, called the Albanese variety of $\HNC/\Gamma$, see \cite[chapter 12]{voisinTheorieHodgeGeometrie2002}. This map $A$ lifts to a quasi-isometric holomorphic map $\widetilde{A}: \HNC/[\Gamma,\Gamma] \to \C^{N}$, which implies that $\HNC/[\Gamma,\Gamma]$ is holomorphically convex. This manifold is compact if and only if $b_1(\Gamma)=0$. There are examples of uniform and torsion-free lattices $\Gamma$ such that $\HNC/[\Gamma,\Gamma]$ is not compact and contains compact subvarieties of positive dimension, see \cite{cartwrightCartwrightStegerSurface2017}. However, for any uniform and torsion-free arithmetic lattice $\Gamma$ of $\PU(n,1)$ with positive first Betti number, there is a finite index subgroup $\Gamma_1$ of $\Gamma$ such that the Albanese map $A_1$ of $\HNC/\Gamma_1$ is an immersion, see \cite[Section 3.1]{llosaisenrichSubgroupsHyperbolicGroups2024} or \cite{eyssidieuxOrbifoldKahlerGroups2018}. In particular, the lift $\widetilde{A_1}:\HNC/[\Gamma_1,\Gamma_1] \to \C^{N_1}$ of $A_1$ is an immersion, where $N_1=\frac{b_1(\Gamma_1)}{2}$, and we deduce that $\HNC/[\Gamma_1,\Gamma_1]$ does not contain a compact subvariety of positive dimension. Thus $\HNC/[\Gamma_1,\Gamma_1]$ is a Stein manifold in these cases.
		\end{enumerate}
		
	\end{examples}
	
	To conclude this section, let us note that the question of finding sufficient conditions on the group $\Gamma$ for the quotient $X/\Gamma$ to be Stein admits a natural generalization to the case where $X$ is a higher rank Hermitian symmetric space and $\Gamma$ is a group acting freely and properly discontinuously by holomorphic isometries on $X$. Here is a family of examples of quotients of the bidisk $\mathbb{D}\times \mathbb{D}$ that are easily proven to be Stein manifolds.
	
	\begin{example}
		Let $\Gamma=\pi_1(S_g)$ be a cocompact lattice in $\mathrm{PSL}_2(\R)$. Define the following action of $\Gamma$ on the bidisk $\mathbb{D}\times \mathbb{D}$\begin{equation*}
			\forall \gamma \in \Gamma, \forall (z,w)\in \mathbb{D}\times \mathbb{D},\ \gamma\cdot (z,w) := (\gamma \cdot z,\overline{\gamma \cdot \overline{w}}).
		\end{equation*}
		Then $\Delta := \{(z,\overline{z}) \ \vert \ z\in \mathbb{D}\}$ is a totally real geodesically convex subset of $\mathbb{D}\times \mathbb{D}$, on which $\Gamma$ acts cocompactly. Using \cite[Proposition 3.2]{chenDiscreteGroupsHolomorphic2013}, we get that $(\mathbb{D}\times \mathbb{D})/\Gamma$ is a Stein manifold.
	\end{example}
	
	\section{Discrete parabolic subgroups of $\PU(n,1)$}\label{section-paraboliques}
	
	This section is organised as follows. We first recall the definition of the complex hyperbolic distance on the ball, and describe the stabilizer of a point at infinity. Then we state and prove Theorem \ref{thm-caracterisation-parab-stein} which characterises the discrete and torsion-free parabolic subgroups of $\PU(n,1)$ for which $\HNC/\Gamma$ is a Stein manifold, and which implies Theorem \ref{thm-parabolique-abelien-delta-stein}-(b). We then show that if $\Gamma$ is a discrete parabolic subgroup which satisfies $\delta(\Gamma)<2$ or preserves a totally real geodesic submanifold of $\HNC$, then $\Gamma$ is virtually Abelian, thus completing the proof of Theorem \ref{thm-parabolique-abelien-delta-stein}. Afterwards we give an example of a discrete parabolic subgroup with $\delta(\Gamma)=2$ and for which $\HNC/\Gamma$ is not Stein. We also construct a complex hyperbolic bundle of punctured disks over a non-compact Cousin manifold. This complex hyperbolic bundle is not holomorphically convex, but holomorphic functions separate points by \cite[Theorem 1.1]{miebachQuotientsBoundedHomogeneous2024}. Notice that a parabolic quotient $\HNC/\Gamma$ is Stein if and only if it is holomorphically convex, as follows from Proposition \ref{prop-criteres}-(a).
	
	\subsection{The parabolic biholomorphisms of the ball}\label{ss-section-preli-parab}
	
	Let $h$ be the Hermitian form on $\C^{n+1}$ associated with the quadratic form\begin{equation*}
		q(z_1,\dots,z_{n+1}):=-\lvert z_1\rvert^{2} + \sum_{i=2}^{n+1}\lvert z_i\rvert^{2},
	\end{equation*}
	and let $[\cdot]:\C^{n+1}\setminus\{0\} \to \PNC$ denote the projection onto the complex projective space $\PNC$. The open subset of $\PNC$ defined by\begin{equation*}
		\HNC := \{[v]\in \PNC \ \vert \ q(v)<0\}
	\end{equation*}
	is biholomorphic to the unit ball of $\C^{n}$. It can be endowed with a complete Kähler metric of negative sectional curvature pinched between $-4$ and $-1$, for which the distance between two points $x,y\in \HNC$ is given by the formula\begin{equation*}
		\cosh^{2}d(x,y) = \frac{h(\widetilde{x},\widetilde{y}) h(\widetilde{y},\widetilde{x})}{ h(\widetilde{x},\widetilde{x}) h(\widetilde{y},\widetilde{y})},
	\end{equation*}
	where $\widetilde{x},\widetilde{y}\in \C^{n+1}$ denote lifts of $x,y$. Moreover, every biholomorphism of the ball is an isometry for this metric, and the group of biholomorphic isometries of $\HNC$ is isomorphic to $\PU(n,1)$. The action of $\PU(n,1)$ on $\HNC$ extends to an action by homeomorphisms on the closed ball $\HNC\cup \PHNC$, where\begin{equation*}
		\PHNC := \{[v]\in \PNC \ \vert \ q(v)=0\}.
	\end{equation*}
	
	Let us fix a point $\xi\in\PHNC$. There exists a basis $f_{\xi}=(f_1,f_2,e_1,\dots,e_{n-1})$ of $\C^{n+1}$ such that $\xi=[f_1]$ and in which the quadratic form $q$ has the following expression\begin{equation*}
		q\left(\alpha f_1+\beta f_2+\sum_{i=1}^{n-1}u_ie_i\right) = 2\Re(\alpha\overline{\beta}) + \sum_{i=1}^{n-1}\lvert u_i\rvert^{2}.
	\end{equation*}
	The biholomorphism 
	\begin{equation*}
		\left\{\begin{array}{rcl}
			\{(\alpha,u)\in \C\times\C^{n-1}\ \vert \ 2\Re(\alpha)+\lVert u\rVert^{2}<0\} &\longrightarrow &\HNC\\
			(\alpha,u) &\longmapsto & [\alpha f_1+f_2+u]
		\end{array}	\right.
	\end{equation*}
	defines a global chart of $\HNC$, in which Busemann functions at $\xi$ are the translates of the function $B_{\xi}$ defined by
	\begin{equation*}
		e^{2B_{\xi}(\alpha,u)}=\frac{-2}{2\Re(\alpha)+\lVert u\rVert^{2}}.
	\end{equation*}
	In the basis $f_{\xi}$, let us define three subgroups $M,A$ and $N$ of $\PU(n,1)$ by the associated groups of matrices \begin{align*}
		\mathcal{M} &= \left\{\left(\begin{array}{ccc}
			1&0&0\\
			0&1&0\\
			0&0&T
		\end{array}\right)\lvert\ T\in \mathrm{U}(n-1)\right\}, \\
		\mathcal{A} &= \left\{\left(\begin{array}{ccc}
			e^{t}&0&0\\
			0&e^{-t}&0\\
			0&0&\mathrm{I}_{n-1}
		\end{array}\right)\lvert\ t\in \R\right\}, \\
		\mathcal{N} &=\left\{	\left(\begin{array}{ccc}
			1 & a & -{}^{t}\overline{b} \\
			0 & 1 & 0 \\
			0 & b & \mathrm{I}_{n-1}
		\end{array}\right)\lvert \ b\in \C^{n-1},a\in \C, \lVert b\rVert^{2} = -2\Re(a)\right\}.
	\end{align*}
	Then the stabilizer of a point $\xi\in \PHNC$ in $\PU(n,1)$ decomposes as\begin{equation*}
		\operatorname{Stab}_{\xi}(\HNC)= MAN.
	\end{equation*}
	For $T\in \mathrm{U}(n-1),b\in \C^{n-1}$ and $c\in \R$, let $(T,b,c)$ denote the element of the group $MN$ defined in the basis $f_{\xi}$ by the matrix\begin{equation*}
		\left(\begin{array}{ccc}
			1 & -\frac{\lVert b\rVert^{2}}{2}+ic & -(\langle Te_j,b\rangle)_{1\le j\le n-1} \\
			0 & 1 & 0 \\
			0 & b & T
		\end{array}\right).
	\end{equation*}
	The group law on $MN$ is given by\begin{equation*}
		(T,b,c)\cdot(T',b',c') = (TT',b+Tb',c+c'+\Im\langle b,Tb' \rangle),
	\end{equation*}
	where $\langle\cdot,\cdot\rangle$ is the standard Hermitian product on $\C^{n-1}$. This group is a semi-direct product $\mathrm{U}(n-1)\ltimes N$, and $N$ is isomorphic to the Heisenberg group of real dimension $2n-1$. The center $Z(N)$ of $N$ is the set of elements of the form $(\Id,0,c)$, with $c\in \R$. We denote by $\pi$ the projection of $\mathrm{U}(n-1)\ltimes N$ onto $\mathrm{U}(n-1)$, and by $\Pi$ the morphism $\mathrm{U}(n-1)\ltimes N \to \mathrm{U}(n-1)\ltimes \C^{n-1}$ which to $(T,b,c)\in \mathrm{U}(n-1)\ltimes N$ associates the holomorphic isometry $z\mapsto Tz+b$ of $\C^{n-1}$.
	
	A subgroup $\Gamma$ of $\PU(n,1)$ acts freely and properly discontinuously on $\HNC$ if and only if it is torsion-free and discrete in $\PU(n,1)$. It is parabolic if it fixes a point $\xi$ in $\PHNC$ and if all the eigenvalues of its elements are of modulus 1, which amounts to saying that, in the model described above, $\Gamma$ is a subgroup of $\mathrm{U}(n-1)\ltimes N$, or equivalently that $\Gamma$ preserves horoballs at $\xi$. We write for all $\gamma\in \Gamma$\begin{align*}
		&\gamma = (\pi(\gamma),b(\gamma),c(\gamma)), \text{ with } \pi(\gamma)\in \mathrm{U}(n-1),b(\gamma)\in \C^{n-1},c(\gamma)\in \R, \text{ and }\\
		&\Pi(\Gamma) = (\pi(\gamma),b(\gamma)).
	\end{align*}
	The parabolic group $\Gamma$ is said to be \textit{unipotent} if $\pi(\Gamma)$ is trivial.
		
	\subsection{A characterization of Stein parabolic quotients of the ball}
	
	As explained above, we identify a parabolic subgroup of $\PU(n,1)$ with a subgroup of $\mathrm{U}(n-1)\ltimes N$, and we denote by $\pi$, respectively $\Pi $, the projection of $\mathrm{U}(n-1)\ltimes N$ onto $\mathrm{U}(n-1)$, respectively $ \mathrm{U}(n-1)\ltimes \C^{n-1}$. The following lemma is probably classical, and we give its proof for the reader's convenience.
	
	\begin{lemma}\label{lem-ssgrp-pi-abelien}
		Let $\Gamma$ be a discrete and torsion-free parabolic subgroup of $\PU(n,1)$. Then there exists a finite-index subgroup $\Gamma_1$ of $\Gamma$ such that $\Pi(\Gamma_1)$ is Abelian.
	\end{lemma}
	
	\begin{proof}[Proof of Lemma \ref{lem-ssgrp-pi-abelien}]
		As in Subsection \ref{ss-section-preli-parab}, let us fix a basis of $\C^{n+1}$ which induces an identification between $\Gamma$ and a discrete subgroup of $\mathrm{U}(n-1)\ltimes N$. By Margulis Lemma, $\Gamma$ is virtually nilpotent. The existence of $\Gamma_1$ follows from the classical fact that a nilpotent subgroup of $\mathrm{U}(n-1)\ltimes \C^{n-1}$ is virtually Abelian. To show this fact, the first observation, that we will not prove here, is that a nilpotent subgroup of $\mathrm{U}(n-1)$ is virtually Abelian. Let $\Gamma_1$ be a finite-index nilpotent subgroup of $\Gamma$ such that $\pi(\Gamma_1)$ is Abelian. Seeking a contradiction, let us assume that the nilpotent group $\Pi(\Gamma_1)$ is not Abelian. There is a non-trivial element $z$ in the center of $\Pi(\Gamma_1)$ which can be written as a product of commutators $z=[x_1,y_1]\dots[x_k,y_k]$, with $x_1,\dots,x_k,y_1,\dots,y_k\in \Pi(\Gamma_1)$. Let us write\begin{align*}
			&x_i=(\pi(x_i),b(x_i)),\\
			&y_i=(\pi(y_i),b(y_i)),\\
			&z= (\Id,b(z)), 
		\end{align*}
		with $\pi(x_i),\pi(y_i)\in \mathrm{U}(n-1)$ which commute, and $b(x_i),b(y_i),b(z)\in \C^{n-1}$. We will now show that $b(z)=0$, which means that $z$ is trivial, a contradiction.
		Given $\pi_1,\pi_2\in \mathrm{U}(n-1)$ which commute and $b_1,b_2\in \C^{n-1}$, we compute that \begin{align*}
			[(\pi_1,b_1),(\pi_2,b_2)]&=(\pi_1,b_1)\ (\pi_2,b_2)\ (\pi_1^{-1},-\pi_1^{-1}b_1)\ (\pi_2^{-1},-\pi_2^{-1}b_2) \\
			&= (\pi_1\pi_2,b_1+\pi_1b_2)(\pi_1^{-1}\pi_2^{-1},-\pi_1^{-1}b_1 - \pi_1^{-1}\pi_2^{-1}b_2)\\
			&= (\Id,b_1+\pi_1b_2 +(\pi_1\pi_2)(-\pi_1^{-1}b_1 - \pi_1^{-1}\pi_2^{-1}b_2))\\
			&= (\Id,(\Id-\pi_2)b_1-(\Id-\pi_1)b_2).
		\end{align*}
		Therefore
		\begin{equation}\label{b31}
			b(z)=\sum_{i=1}^{k}(\Id-\pi(y_i))b(x_i)-(\Id-\pi(x_i))b(y_i).
		\end{equation}
		Moreover, $z$ commutes with all elements of $E:=\{x_1,\dots,x_k,y_1,\dots,y_k\}$, which implies that
		\begin{equation}\label{b32}
			\forall \gamma \in E,\quad \pi(\gamma)b(z)=b(z).
		\end{equation}
		Choose a basis $e=(e_1,\dots,e_{n-1})$ that diagonalizes all elements of $\pi(E)$ and express $\pi(\gamma)$ in this basis as $\operatorname{Diag}(a_1(\gamma),\dots,a_{n-1}(\gamma))$ for $\gamma\in E$. For all $j\in\{1,\dots,n-1\}$, if there exists $\gamma\in E$ such that $a_j(\gamma)\ne 1$, then the $j^{\text{th}}$ coordinate $b_j(z)$ of $b(z)$ in the basis $e$ must vanish according to Formula \eqref{b32}, and if $a_j(\gamma)= 1$ for all $\gamma\in E$, then $b_j(z)=0$ according to Formula \eqref{b31}. Thus $b(z)=0$, which gives the contradiction we were looking for and proves that $\Pi(\Gamma_1)$ is Abelian.
	\end{proof}
	
	\begin{remark}
		It is a classical fact that any discrete and torsion-free parabolic subgroup of $\PO(n,1)$ is virtually Abelian. This can be seen using Lemma \ref{lem-ssgrp-pi-abelien}, because for such a group, we have $c(\gamma)=0$ and thus $\Pi(\Gamma)$ is isomorphic to $\Gamma$.
	\end{remark}
	
	Let $\Gamma$ and $\Gamma_1$ be as in Lemma \ref{lem-ssgrp-pi-abelien}. Set\begin{equation*}
		V_1:= \bigcap_{\gamma\in \Gamma_1}\ker(\operatorname{Id}-\pi(\gamma)),
	\end{equation*} and let $p:\C^{n-1}\to V_1$ be the orthogonal projection onto $V_1$. Finally, define \begin{equation*}
	W_1 := \operatorname{Span}(\{p(b(\gamma)) \ \vert \ \gamma\in \Gamma_1\}),
	\end{equation*}
	where $b(\gamma)=\Pi(\gamma)\cdot 0 \in \C^{n-1}$. In the following statement, a linear subspace $W$ of $\C^{n-1}$ is said to be totally real if $W\cap iW=\{0\}$.
	
	\begin{theorem}\label{thm-caracterisation-parab-stein}
		Let $\Gamma$ be a discrete and torsion-free parabolic subgroup of $\PU(n,1)$, and let $\Gamma_1,p, V_1$ and $W_1$ be as above. Then $\HNC/\Gamma$ is a Stein manifold if and only if $W_1$ is totally real.
	\end{theorem}
		
	\begin{proof}[Proof of Theorem \ref{thm-caracterisation-parab-stein}] Let $\Gamma_1$ be as in Lemma \ref{lem-ssgrp-pi-abelien}. Since $\HNC/\Gamma_1$ is a finite covering of $\HNC/\Gamma$, one of these two manifolds is Stein if and only if the other one is. To simplify the notation, we can thus assume without loss of generality that $\Gamma_1=\Gamma$.
		
		\underline{Step 1.} There is a separation of $\Gamma$ into an elliptic and a unipotent part.
		
		 The group $\pi(\Gamma)$ is Abelian, hence there is an orthonormal basis $e=(e_1,\dots,e_{n-1})$ of $\C^{n-1}$ as well as morphisms $a_1,\dots,a_{n-1}$ from $\Gamma$ to the unit circle in $\C$ such that, in the basis $e$
		\begin{equation*}
			\pi(\gamma) = \operatorname{Diag}(a_1(\gamma),\dots,a_{n-1}(\gamma)).
		\end{equation*}
		For $\gamma \in \Gamma$, let $(b_1(\gamma),\dots,b_{n-1}(\gamma))$ be the coordinates of $b(\gamma)$ in the basis $e$.
		Up to permuting the elements of $e$, we can assume that $(e_{k+1},\dots,e_{n-1})$ forms a basis of $V_1$ for some integer $k\in \{0,\dots,n-1\}$. For all $i\in \{1,\dots,k\}$, there is an element $\gamma_i \in \Gamma$ such that $a_i(\gamma_i)\ne 1$. Set
		\begin{equation*}
			\lambda_i := \frac{b_i(\gamma_i)}{1-a_i(\gamma_i)}.
		\end{equation*}
		As $\Pi(\Gamma)$ is Abelian, we have
		\begin{equation*}
			\forall \gamma \in \Gamma, \forall i\le k, b_i(\gamma)=\lambda_i(1-a_i(\gamma)).
		\end{equation*}
		We set $\Lambda := {}^{t}(\lambda_1,\dots,\lambda_k,0,\dots,0) \in \C^{n-1}$, so that
		\begin{equation*}
			\forall \gamma \in \Gamma, (\pi(\gamma)-\Id)\Lambda = {}^{t}(-b_1(\gamma),\dots,-b_k(\gamma),0,\dots 0).
		\end{equation*}
		Let $T_{\Lambda}$ be the element $(\Id,\Lambda,0)$ of $\mathrm{U}(n-1)\ltimes N$ and, for all $\gamma \in \Gamma$, define $\phi(\gamma) :=T_{\Lambda}^{-1}\gamma T_{\Lambda}$. A computation shows that for all $\gamma\in \Gamma$, we have
		\begin{equation*}
			\phi(\gamma) = (\pi(\gamma),b(\phi(\gamma)),c(\phi(\gamma))),
		\end{equation*}
		where $b(\phi(\gamma)) := b(\gamma) + (\pi(\gamma)-\Id)\Lambda \in V_1$ and $c(\phi(\gamma)) \in \R$. Set
		\begin{equation*}
			\phi(\gamma)_e := (\pi(\gamma),0,0) \text{ and }
			\phi(\gamma)_u := (\Id,b(\phi(\gamma)),c(\phi(\gamma))).
		\end{equation*}
		Using that $b(\phi(\gamma))\in V_1$ for all $\gamma\in \Gamma$, we get
		\begin{equation*}
			\forall \gamma,\gamma'\in \Gamma, [\phi(\gamma)_e,\phi(\gamma')_u] = \Id.
		\end{equation*}It is also easily verified that $\phi(\gamma) = \phi(\gamma)_e\phi(\gamma)_u$ for all $\gamma \in \Gamma$. 
		More generally, $U(k)$, seen as a subgroup of $\mathrm{U}(n-1)$ fixing $V_1$ pointwise, commutes with $\phi(\gamma)_u$ for $\gamma\in \Gamma$. In particular, $\phi(\Gamma)_E:= \{\phi(\gamma)_e \ \vert \ \gamma \in \Gamma\}$ and $\phi(\Gamma)_U := \{\phi(\gamma)_u \ \vert \ \gamma \in \Gamma\}$ are groups and $\phi(\Gamma)_U<N$. Moreover, $\phi(\Gamma)_U$ is discrete in $N$. Indeed, let $(\gamma_k)_{k\in \N}$ be a sequence in $\Gamma$ such that $\phi(\gamma_k)_u \underset{k\to +\infty}{\longrightarrow} \Id$. After passing to a subgroup, we can assume that the sequence $\phi(\gamma_k)_e$ converges to a limit $M\in \mathrm{U}(n-1)$. Thus $\phi(\gamma_k)$ converges to $M$, and since $\phi(\Gamma)$ is discrete, this sequence has to be stationary. Hence $\phi(\Gamma)_U$ is discrete in $N$.
		
		\underline{Step 2} We deduce the characterization of Stein quotients.
		
		We can rewrite $W_1$ as
		\begin{equation*}
			W_1 = \operatorname{Span}_{\R}(\{b(\phi(\gamma)) \ \vert \ \gamma \in \Gamma\}).
		\end{equation*}
		From \cite[Theorem 1.4]{miebachQuotientsBoundedHomogeneous2024}, we obtain that the quotient $\HNC/\phi(\Gamma)_U$ is Stein if and only if $W_1$ is totally real.
		
		Assume that $W_1$ is totally real. Then $\HNC/\phi(\Gamma)_U$ is a Stein manifold, so it has a strictly plurisubharmonic exhaustion function $\psi_U:\HNC/\phi(\Gamma)_U\to \R_+$. Moreover, the holomorphic action of $U(k)$ on $\HNC$ descends to the quotient $\HNC/\phi(\Gamma)_U$, and by averaging $\psi_U$ over the orbits of $U(k)$, we can assume that $\psi_U$ is $U(k)$-invariant. Then $\psi_U$ lifts to a strictly plurisubharmonic function $\widetilde{\psi_U}:\HNC\to \R_+$ which is invariant by $\phi(\Gamma)_U$ and $U(k)$. This function descends to a strictly plurisubharmonic function $\psi:\HNC/\phi(\Gamma)\to \R_+$, which is an exhaustion. Thus, $\HNC/\phi(\Gamma)$, and therefore $\HNC/\Gamma$, are Stein.
		
		Conversely, if $\HNC/\Gamma$, hence $\HNC/\phi(\Gamma)$, is a Stein manifold, then $\HNC/\phi(\Gamma)$ admits a strictly plurisubharmonic exhaustion function. This implies that $\HNC/\phi(\Gamma)_U$ has a strictly plurisubharmonic exhaustion function, thus $\HNC/\phi(\Gamma)_U$ is a Stein manifold. Therefore, $W_1$ is totally real.
	\end{proof}
	
	Using Theorem \ref{thm-caracterisation-parab-stein}, we obtain as a corollary the second point of Theorem \ref{thm-parabolique-abelien-delta-stein}.
	
	\begin{corollary}\label{cor-virt-ab-stein}
		Let $\Gamma$ be a discrete and torsion-free parabolic subgroup of $\PU(n,1)$. If $\Gamma$ is virtually Abelian, then $\HNC/\Gamma$ is a Stein manifold.
	\end{corollary}
	
	\begin{proof}
		As in the proof of Theorem \ref{thm-caracterisation-parab-stein}, we assume without loss of generality that $\Gamma$ is Abelian and we decompose $\phi(\Gamma)$ into an elliptic part $\phi(\Gamma)_E$ and a unipotent part $\phi(\Gamma)_U$. Then $\phi(\Gamma)_U$ is a discrete and Abelian parabolic subgroup of $\PU(n,1)$. It is known that the quotient of the complex hyperbolic space by such a subgroup is a Stein manifold, see \cite{chenDiscreteGroupsHolomorphic2013}. This is also a particular case of \cite[Theorem 1.4]{miebachQuotientsBoundedHomogeneous2024}, because with the notations of Step 1 above, it can be verified that for all $\gamma,\gamma'\in \Gamma$, the identity $[\phi(\gamma),\phi(\gamma')]=\Id$ implies that
		\begin{equation*}
			\Im\langle b(\phi(\gamma)),b(\phi(\gamma'))\rangle = 0,
		\end{equation*}
		and this implies that $W_1$ is totally real. Hence $\HNC/\Gamma$ is a Stein manifold.
	\end{proof}
	\subsection{Proof of Theorem \ref{thm-parabolique-abelien-delta-stein}}
	
	We first recall a formula for the critical exponent of a discrete and torsion-free parabolic subgroup $\Gamma$ of $\PU(n,1)$, for which we refer to \cite{corletteLimitSetsDiscrete1999} or \cite[§3]{dalboSeriesPoincareGroupes2000}. Let $\Gamma_1$ be a finite-index subgroup of $\Gamma$ such that $\Pi(\Gamma_1)$ is Abelian. Define $l\in\{0,1\}$ as the dimension of the real subspace spanned by $Z(N)\cap \Gamma_1$, where $Z(N)\simeq\R$ is the center of $\R$, and $k\in\{0,\dots 2n-2\}$ as the dimension of the subspace of $\C^{n-1}$ spanned by $\{b(\gamma) \mid \gamma \in \Gamma_1\}$. Then\begin{equation}\label{eq-delta-parab}
		\delta(\Gamma) = \frac{2l+k}{2},
	\end{equation}
	see \cite[Proof of Lemma 3.5]{corletteLimitSetsDiscrete1999}
	
	\begin{proof}[Proof of Theorem \ref{thm-parabolique-abelien-delta-stein}]
		The second point of the theorem is given by Corollary \ref{cor-virt-ab-stein}. For the first point, let $\Gamma$ be a discrete and torsion-free parabolic subgroup of $\PU(n,1)$ which is not virtually Abelian. We will show that $\delta(\Gamma)\ge2$ by finding two elements $x,y\in \Gamma$ that generate a group of critical exponent equal to 2. Let us fix, as in Subsection \ref{ss-section-preli-parab}, a basis $f$ of $\C^{n+1}$ which induces an identification between $\Gamma$ and a subgroup of $\mathrm{U}(n-1)\ltimes N$. Let $\Gamma_1$ be a finite-index subgroup of $\Gamma$ such that $\Pi(\Gamma_1)$ is Abelian, given by Lemma \ref{lem-ssgrp-pi-abelien}. Since the set of commutators of elements of $\Gamma_1$ is included in the kernel of $\Pi$, which coincides with the center $Z(N)$ of $N$, and $\Gamma_1$ is not Abelian, we deduce that $\Gamma$ contains two elements $x$ and $y$ such that $\Pi(x)$ and $\Pi(y)$ commute, but $x$ and $y$ do not. Then according to Formula \eqref{eq-delta-parab}, the critical exponent of the group generated by $x$ and $y$ is $\frac{2l+k}{2}$, where $l\in\{0,1\}$ is the dimension of the $\R$-span of the elements $c(\gamma)$ for $\gamma\in\Gamma$ and $k\in\{0,1,2\}$ is the dimension of the $\R$-span of $b(x)$ and $b(y)$. Since $x$ and $y$ do not commute, we see that $l=1$ and $k=2$. Thus $\delta(\langle x,y\rangle)=2$.
		
		Now assume that $\Gamma$ preserves a a totally real and totally geodesic submanifold of $\HNC$. Then we can realize $\Gamma$ as a discrete and virtually nilpotent subgroup of \begin{equation*}
			\mathbb{P}((O(k-1)\ltimes \R^{k-1})\times U(n-k))
		\end{equation*}for some integer $k\in\{1,\dots,n\}$. Consequently, $\Gamma$ is virtually Abelian (see the remark after the proof of Lemma \ref{lem-ssgrp-pi-abelien}). We deduce from Theorem \ref{thm-caracterisation-parab-stein} that $\HNC/\Gamma$ is a Stein manifold.
	\end{proof}
	
	\begin{corollary}\label{cor-parab-dim2}
		Let $\Gamma$ be a discrete and torsion-free parabolic subgroup of $\PU(2,1)$. Then $\mathbb{H}^{2}_{\C}/\Gamma$ is Stein if and only if $\Gamma$ is virtually Abelian.
	\end{corollary}
	
	\begin{proof}
		If $\Gamma$ is not virtually Abelian, choose two elements $x,y\in \Gamma$ as in the proof of Theorem \ref{thm-parabolique-abelien-delta-stein}. Since $x$ and $y$ do not commute, we get that $\Im\langle b(x),b(y)\rangle \ne 0$, and thus $W:=\operatorname{Span}_{\R}(b(x),b(y))\subset\C$ is equal to $\C$. In particular $W$ is not totally real and using Theorem \ref{thm-caracterisation-parab-stein}, we deduce that $\HNC/\Gamma$ admits a covering $\HNC/\langle x_0,y_0\rangle$ which is not Stein. As any covering of a Stein manifold is Stein, see \cite{steinUberlagerungenHolomorphvollstandigerKomplexer1956}, this implies that $\HNC/\Gamma$ is not Stein.
	\end{proof}
	
	\subsection{Examples of parabolic quotients of the ball}
	In the following two examples, we fix a basis $f=(f_1,f_2,e_1,\dots,e_{n-1})$ of $\C^{n+1}$ as in Subsection \ref{ss-section-preli-parab}, which induces an identification between parabolic subgroups of $\PU(n,1)$ fixing $[f_1]$ and subgroups of $\mathrm{U}(n-1)\ltimes N$. We recall from \cite[Proposition 1.1]{miebachQuotientsBoundedHomogeneous2024} that if $\Gamma$ is a unipotent discrete and torsion-free subgroup of $\PU(n,1)$, then $\HNC/\Gamma$ is holomorphically separable.
	
	\begin{example}Here is an example of a discrete unipotent subgroup $\Gamma$ of $\PU(n,1)$ with $\delta(\Gamma)=2$, for which $\HNC/\Gamma$ is not holomorphically convex. The group $\Gamma$ generated by $\gamma_1 :=(\Id,e_1,0)$ and $\gamma_2:=(\Id, ie_1,0)$ is the set of all elements of the form \begin{equation*}
			(\Id,(k_1+ik_2)e_1,2\ell-k_1k_2),
		\end{equation*}
		where $(k_1,k_2,\ell)\in \Z^3$. In particular, $\Gamma$ is discrete, and Formula \eqref{eq-delta-parab} shows that $\delta(\Gamma)=2$. Finally, the quotient $\HNC/\Gamma$ naturally identifies with a bundle of punctured disks over the base $B:=\C/(\Z+i\Z) \times \C^{n-2}$. If $\HNC/\Gamma$ were holomorphically convex, it would be Stein by Proposition \ref{prop-criteres}-(a) and we would deduce that $B$ is a Stein manifold by \cite[Lemma 1.6]{coltoiuOpenSetsStein1997}, which is not the case. Therefore, $\HNC/\Gamma$ is not holomorphically convex.
	\end{example}
	
	\begin{example}Here is an example of a complex hyperbolic bundle of punctured disks over a Cousin group. We work in dimension $n=3$, but this example generalizes to any dimension $n\ge3$. Using the identification introduced before the previous example, define three vectors in $\C^{2}=\C e_1\oplus \C e_2$ by $v_1=e_1$, $v_2=e_2$, and $v_3=ae_1+be_2$ for $(a,b)\in \C^{2}$ two complex numbers such that\begin{equation*}
			\left\{\begin{array}{l}
				\lambda:=\Im(a)=\Im(b)\ne 0, \\
				\Re(a)-\Re(b) \notin \mathbb{Q}.
			\end{array}\right.
		\end{equation*}
		The fact that $\Im(a)\ne0$ and $\Im(b)\ne0$ implies that $v_1,v_2$ and $v_3$ are $\R$-linearly independent, and both conditions together imply that $1,a$ and $b$ are $\Z$-linearly independent. We deduce that the subgroup $\Gamma_0$ of $\C^{2}$ generated by $v_1,v_2$ and $v_3$ is discrete, and that the quotient $\C^{2}/\Gamma_0$ has no compact factor and does not admit any non-constant holomorphic function (see for example \cite[pages 451-452]{napierConvexityPropertiesCoverings1990}).
		Let $\Gamma$ be the subgroup of $\mathrm{U}(n-1)\ltimes N$ generated by the three elements $\gamma_i=(\Id,v_i,0)$ for $i=1,2$ and $3$.
		The equality $\Im(a)=\Im(b)=\lambda$ implies that
		\begin{equation*}
			[\gamma_3,\gamma_1]=[\gamma_3,\gamma_2] = (\Id,0,2\lambda).
		\end{equation*}
		Any element of $\Gamma$ is of the form $\gamma_1^{k_2}\gamma_2^{k_2}\gamma_3^{k_3}[\gamma_3,\gamma_1]^{\ell}$,
		with $(k_1,k_2,k_3,\ell)\in \Z^{4}$, and we deduce that $\Gamma$ is the set of all elements of the form\begin{equation*}
			(\Id,k_1v_2+k_2v_2+k_3v_3,((k_1+k_2)k_3+2\ell)\lambda).
		\end{equation*}
		with $(k_1,k_2,k_3,\ell)\in \Z^{4}$. Consequently, $\Gamma$ is discrete, and $\mathbb{H}^{3}_{\C}/\Gamma$ is biholomorphic to a bundle of punctured disks over $\C^{2}/\Gamma_0$. Since $\C^{2}/\Gamma_0$ is not Stein, we deduce as in the previous example that $\mathbb{H}^{3}_{\C}/\Gamma$ is not holomorphically convex. Additionally, Formula \eqref{eq-delta-parab} shows that $\delta(\Gamma)=\frac{5}{2}$.
	\end{example}
	
	\section{Holomorphic convexity and geometrically finite \\subgroups}\label{section-hol-cvx-geom-fini}
	
	In this section, we prove Theorems \ref{thm-quotient-hol-cvx} and \ref{thm-geom-fini-delta-stein}.
	
	\subsection{Proof of Theorem \ref{thm-quotient-hol-cvx}}
	
	We will need the following lemma about parabolic quotients of the ball.
	
	\begin{lemma}\label{lem-quotient-horoboule-stein}
		Let $P$ be a discrete and torsion-free parabolic subgroup of $\PU(n,1)$, and $\xi\in \PHNC$ the point fixed by $P$. The following statements are equivalent:
		\begin{enumerate}
			\item $\HNC/P$ is a Stein manifold.
			\item For any horoball $H\subset \HNC$ at $\xi$, the quotient $H/P$ is a Stein manifold.
			\item There exists a horoball $H\subset \HNC$ at $\xi$ for which $H/P$ is a Stein manifold.
		\end{enumerate}
	\end{lemma}
	
	\begin{proof}
		We first show the implication 1 $\implies$ 2. If $\HNC/P$ is a Stein manifold, it admits a strictly plurisubharmonic exhaustion function $\psi:\HNC/P\to \R_+$. The Busemann function  $B:\HNC\to \R$ at $\xi$ is invariant under the action of $P$. Let $H_{\lambda} := B^{-1}((-\infty,\lambda))$ be a horoball at $\xi$. The function $\psi + \frac{1}{\lambda - B}$ defined on $H_{\lambda}/P$ is a strictly plurisubharmonic exhaustion function on $H_{\lambda}/P$, which shows that $H_{\lambda}/P$ is a Stein manifold.
		
		The implication 2 $\implies 1$ is a direct consequence of the last theorem of Subsection \ref{ss-section-stein}.
		
		The implication 2 $\implies$ 3 is immediate. We now show that 3 $\implies$ 2. We fix, as in Subsection \ref{ss-section-preli-parab}, a basis $f=(f_1,f_2,e_3,\dots,e_{n+1})$ of $\C^{n+1}$ which induces an identification between parabolic elements of $\PU(n,1)$ fixing $\xi=[f_1]$ and elements of $\mathrm{U}(n-1)\ltimes N$. The elements of $P$, seen as biholomorphisms of $\PNC$, commute with the biholomorphisms $L_t:\PNC\to \PNC$ defined for all real numbers $t$ in the basis $f$ by the matrices
		\begin{equation*}
			\left(\begin{matrix}
				1 &t &0\\
				0 &1 &0\\
				0 &0 &\mathrm{I}_{n-1}
			\end{matrix}\right).
		\end{equation*}
		With the notations of Subsection \ref{ss-section-preli-parab}, it is easily checked that for any pair $(\lambda,\mu)$ of real numbers, the map $L_t$ with $t=e^{-2\lambda}-e^{-2\mu}$ sends the horoball $H_{\lambda}:=B^{-1}((-\infty,\lambda))$ to the horoball $H_{\mu}:=B^{-1}((-\infty,\mu))$. We deduce that the quotients of horoballs at $\xi$ by $P$ are all biholomorphic.
	\end{proof}
	
	We now come to the proof of Theorem \ref{thm-quotient-hol-cvx}.
	
	\begin{proof}[Proof of Theorem \ref{thm-quotient-hol-cvx}]
		(1 $\implies$ 2) Suppose that $\HNC/\Gamma$ admits a plurisubharmonic exhaustion function $\phi:\HNC/\Gamma\to \R$. Let $P$ be a maximal parabolic subgroup of $\Gamma$. There exists a Busemann function $B$, invariant under $P$, such that the set $C:=B^{-1}((-\infty,0))/P$ is biholomorphic to an open subset of $\HNC/\Gamma$. The function $\phi\vert_{ C} + \frac{-1}{B}$ is a strictly plurisubharmonic exhaustion function of $C$, and therefore $C$ is a Stein manifold. Using Lemma \ref{lem-quotient-horoboule-stein}, we deduce that $\HNC/P$ is a Stein manifold. If now $P$ is any parabolic subgroup of $\Gamma$, it is contained in a maximal parabolic subgroup $P_0$ of $\Gamma$. The manifold $\HNC/P$ is a covering of $\HNC/P_0$, which is Stein, and therefore $\HNC/P$ is Stein.

		(2 $\implies$ 3) This proof is inspired by \cite[Proof of Theorem 1.4]{chenDiscreteGroupsHolomorphic2013}. Suppose that, for every maximal parabolic subgroup $P<\Gamma$, the quotient $\HNC/P$ is a Stein manifold. Recall from Subsection \ref{sec:gen-geom-finite} that the manifold $X_{\Gamma}:=\HNC/\Gamma$ decomposes as
		\begin{equation*}
			X_{\Gamma} =: Q \cup \bigcup_{i=1}^{k}E_i,
		\end{equation*}
		where $Q$ is relatively compact in $X_{\Gamma} \cup \partial X_{\Gamma}$, $k$ is an integer, and for $i\in \{1,\dots,k\}$, each $E_i$ is an open subset in $X_{\Gamma}$ biholomorphic to the quotient of a horoball $B_i^{-1}((-\infty,0))$ by a maximal parabolic subgroup $P_i$ of $\Gamma$, for some Busemann function $B_i$. We also define $E_i'\subset E_i$ as the quotient of the horoball $B_i^{-1}((-\infty,-1))$ by $P_i$. Moreover, the convex core $C_{\Gamma}$ of $X_{\Gamma}$ has compact intersection with $Q$. By the arguments given in the proofs of Proposition \ref{prop-KH-ss-ens} and Lemma \ref{lem-cvx-psh}, we get that the squared distance function to the convex hull $C(\Gamma)$ of the limit set descends to a convex function $\phi$ on $X_{\Gamma}$, which is strictly convex outside $C_{\Gamma}$. By Richberg's theorem, there exists a continuous plurisubharmonic function $\widetilde{\phi}$ which is smooth and strictly plurisubharmonic outside $C_{\Gamma}$, and such that\begin{equation*}
			\phi \le \widetilde{\phi} \le \phi + \frac{1}{2},
		\end{equation*}see \cite[Theorem I.5.21]{demaillyComplexAnalyticDifferential-manuel}. Moreover, Lemma \ref{lem-quotient-horoboule-stein} implies that for any $i\in \{1,\dots,k\}$, the open subset $E_i$ of $X_{\Gamma}$ is a Stein manifold, and admits a strictly plurisubharmonic exhaustion function. Let $\Psi_i$ be a smooth non-negative function that coincides with this function on $E_i'$ and vanishes outside $E_i$. For any integer $j\in \N$, let $T^{i}_j$ be the compact subset of $X_{\Gamma}$ defined by
		\begin{equation*}
			T^{i}_j := \{x\in \overline{E_i}\setminus E'_i \ \vert \ j\le \widetilde{\phi}(x) \le j+1\}.
		\end{equation*}
		Then, as soon as $j\ge1$, the function $\widetilde{\phi}$ is strictly plurisubharmonic on $T_i^{j}$, so there exists a constant $\beta^{i}_j>0$ such that $i\partial\bar{\partial} \Psi_i \ge -\beta^{i}_j i\partial\bar{\partial} \widetilde{\phi}$ on $T^{i}_j$. It follows that there exists a strictly increasing convex function $\lambda:\R_{+}\to\R_{+}$ such that $\lambda(t)\underset{t\to+\infty}{\longrightarrow}+\infty$ and such that
		\begin{equation*}
			N:= \lambda\circ\widetilde{\phi} + \sum_{i=1}^{k}\Psi_i
		\end{equation*}
		is strictly plurisubharmonic on the set
		\begin{equation*}
			\bigcup_{i=1}^{k}\bigcup_{j\ge1}T^{i}_{j}.
		\end{equation*}
		On $Q$, this function $N$ coincides with $\lambda\circ\widetilde{\phi}$ and it is strictly plurisubharmonic on $Q\cap X_{\Gamma}\setminus C_{\Gamma}$. On each $E_i'$, since $\Psi_i$ is strictly plurisubharmonic and $\widetilde{\phi}$ is plurisubharmonic, $N$ is strictly plurisubharmonic. In conclusion, $N$ is strictly plurisubharmonic outside the compact set
		\begin{equation*}
			(C_{\Gamma}\cap Q)\cup\bigcup_{i=1}^{k}T^{i}_{0}.
		\end{equation*}
		Moreover, $N$ is an exhaustion function. Indeed, if $(x_n)_{n\in \N}$ is a sequence in $X_{\Gamma}$ without accumulation point, then, after extracting a subsequence, it converges to the boundary $\partial X_{\Gamma}$ or has values in one of the open sets $E_i'$. In the first case where $x_n \longrightarrow x_{\infty}\in \Omega(\Gamma)/\Gamma$, we claim that $\phi(x_n)\longrightarrow +\infty$. Assuming the contrary, we obtain a sequence $(\widetilde{x_n})_{n\in \N}$ in $\HNC$ converging to an element $\widetilde{x_{\infty}}\in \Omega(\Gamma)$, which remains at bounded distance from $C(\Gamma)$, and thus another sequence $(c_n)_{n\in \N}$ in $C(\Gamma)$ converging to $\widetilde{x_{\infty}}$. Thus, $\widetilde{x_{\infty}}\in \Omega(\Gamma) \cap \partial C(\Gamma)$. This is a contradiction, because $\partial C(\Gamma)=\Lambda(\Gamma)$, see \cite{andersonDirichletProblemInfinity1983}. In the case where the sequence lies in $E_i'$, it does not accumulate and therefore, after passing to a subsequence, we have $\Psi(x_n)\to +\infty$. Thus $N$ is an exhaustion function and therefore $X$ is holomorphically convex.
		
		The implication 3 $\implies$ 1 is classical, see for example \cite[Theorem I.6.14]{demaillyComplexAnalyticDifferential-manuel}.
	\end{proof}
	
	\begin{remark}
		If we replace $\HNC$ by a simply connected complete Kähler manifold $(X,\omega)$ with negatively pinched sectional curvature, and assume that $\Gamma$ is a group acting freely and geometrically finitely by holomorphic isometries on $X$, I do not know if Lemma \ref{lem-quotient-horoboule-stein} remains true (the proof uses the holomorphic maps $L_t$ whose existence is specific to the complex hyperbolic case). In Theorem \ref{thm-quotient-hol-cvx}, it remains true that 1 $\iff$ 3. To show that $1 \implies 3$, one argues as in the proof above, noticing that if $X/\Gamma$ admits a plurisubharmonic exhaustion $\phi$, then the open sets $E_i$ appearing in the decomposition
		\begin{equation*}
			X/\Gamma = Q\cup\bigcup_{i=1}^{k}E_i
		\end{equation*}
		are Stein manifolds, with a strictly plurisubharmonic exhaustion given by $\phi+\frac{-1}{B_i}$, where $B_i$ is a Busemann function on $X$ associated to a parabolic point corresponding to the cusp $E_i$.
	\end{remark}
	
	\subsection{Proof of Theorem \ref{thm-geom-fini-delta-stein}}
	
	For the proof of Theorem \ref{thm-geom-fini-delta-stein}-(a), we will need the following lemma, which is presumably classical, and the proof of which we include for completeness.
	
	\begin{lemma}\label{lem_grp_parab_cyclique_ou_Z2}
		Let $X$ be a complete simply connected Riemannian manifold with negatively pinched curvature, and $P$ a discrete and torsion-free parabolic subgroup of isometries of $X$. Then $P$ is cyclic or contains a copy of $\Z^{2}$.
	\end{lemma}
	
	\begin{proof}
		By Margulis' lemma, $P$ contains a finite-index nilpotent subgroup $P'$. Moreover, $P'$ is finitely generated according to \cite{bowditchDiscreteParabolicGroups1993}. If $P'$ is Abelian, then $P'$ is cyclic or contains a copy of $\Z^{2}$. Since a virtually cyclic torsion-free group is cyclic, we deduce that $P$ is cyclic or contains a copy of $\Z^{2}$. Otherwise, let $g$ be a non-trivial element in the center of $P'$, and $h$ an element of $P'$ which does not belong to the center of $P'$. Then $g$ and $h$ generate a subgroup isomorphic to $\Z$ or $\Z^{2}$. Suppose, by contradiction, that this group is cyclic. Then $g$ and $h$ are powers of an element $k\in P'$. In a torsion-free and finitely generated nilpotent group, the centralizers of an element and its powers coincide, and consequently $g$ and $h$ have the same centralizer in $P'$. This  yields a contradiction, and consequently $P$ contains a copy of $\Z^{2}$.\qedhere
	\end{proof}
	
	\begin{proof}[Proof of Theorem \ref{thm-geom-fini-delta-stein}]
		Let $\Gamma$ be a geometrically finite and torsion-free subgroup of $\PU(n,1)$. Assume first that $\Gamma$ is Gromov-hyperbolic. Then $\Gamma$ does not contain a copy of $\Z^{2}$, see \cite[Corollary III.$\Gamma$.3.10]{bridsonMetricSpacesNonpositive1999}, so according to Lemma \ref{lem_grp_parab_cyclique_ou_Z2}, the non-trivial parabolic subgroups of $\Gamma$ are cyclic. The quotient of the complex hyperbolic space by the action of a cyclic parabolic group is a Stein manifold, as follows, for example, from Theorem \ref{thm-parabolique-abelien-delta-stein}, see also \cite{defabritiisQuotientsUnitBall2001manuel} or \cite{miebachQuotientsBoundedHomogeneous2010}. Theorem \ref{thm-quotient-hol-cvx} implies that $\HNC/\Gamma$ is holomorphically convex.
		
		Suppose now that $\delta(\Gamma)<2$. For any parabolic subgroup $P<\Gamma$, we have $\delta(P)\le \delta(\Gamma)<2$, so $\HNC/P$ is Stein according to Theorem \ref{thm-parabolique-abelien-delta-stein}. Using Theorem \ref{thm-quotient-hol-cvx}, we deduce that $\HNC/\Gamma$ is holomorphically convex. Since $\delta(\Gamma)<2$, this manifold does not contain any compact analytic subvariety of positive dimension according to \cite[Theorem 15]{deyNoteComplexhyperbolicKleinian2020} or Proposition \ref{prop-KH-ss-ens}. We deduce that $\HNC/\Gamma$ is Stein.
		
		Finally, suppose that $\Gamma$ preserves a totally real and totally geodesic submanifold of $\HNC$. Then according to Theorem \ref{thm-parabolique-abelien-delta-stein}, for any parabolic subgroup $P$ of $\Gamma$, the quotient $\HNC/P$ is Stein. Using Theorem \ref{thm-quotient-hol-cvx}, we deduce that $\HNC/\Gamma$ is holomorphically convex. This manifold does not contain any compact analytic subvariety of positive dimension, according for example to Proposition \ref{prop-criteres}-(c). Thus $\HNC/\Gamma$ is Stein.
	\end{proof}

\begin{proof}[Proof of Corollary \ref{cor-free}]
	If $\Gamma$ is free and geometrically finite, then $\HNC/\Gamma$ is holomorphically convex by Theorem \ref{thm-geom-fini-delta-stein}. Using Proposition \ref{prop-criteres}-(e), we deduce that $\HNC/\Gamma$ is Stein.
\end{proof}

	We conclude this section by a remark which follows from Lemma \ref{lem_grp_parab_cyclique_ou_Z2}.
	\begin{remark}
		Let $X$ be a complete simply connected Riemannian manifold with negatively pinched curvature, and $\Gamma$ a group containing no copy of $\Z^{2}$. If $\Gamma$ acts faithfully, discretely, and geometrically finitely by isometries on $X$, then $\Gamma$ is Gromov-hyperbolic. Indeed, a geometrically finite group is hyperbolic relative to its parabolic subgroups. Under the assumption that $\Gamma$ contains no copies of $\Z^{2}$, the parabolic subgroups of $\Gamma$ are cyclic, and in particular Gromov-hyperbolic. This implies that $\Gamma$ itself is Gromov-hyperbolic, see \cite{osinRelativelyHyperbolicGroups2006}. We can use this fact to exhibit finitely generated groups which admit a discrete and faithful representation in $\PU(n,1)$ but no discrete, faithful, and geometrically finite representation of $N$ in $\PU(n,1)$. To do this, let first $\Gamma_0$ be a cocompact arithmetic lattice of the simplest type of $\PU(n,1)$, for the definition of which we refer to \cite[§VIII.5]{borelContinuousCohomologyDiscrete2000}. Then there exists a finite-index torsion-free subgroup $\Gamma<\Gamma_0$ and a morphism $\phi:\Gamma \to \Z$ such that $N:=\ker(\phi)$ is finitely generated but not hyperbolic, see \cite{llosaisenrichSubgroupsHyperbolicGroups2024}, and also \cite{italianoHyperbolic5manifoldsThat2023-manuel} for related results. As a subgroup of $\Gamma$, the group $N$ cannot contain a copy of $\Z^{2}$. Thus, there exists by construction a discrete and faithful representation of $N$ in $\PU(n,1)$, but there is no discrete, faithful, and geometrically finite representation of $N$ in $\PU(n,1)$.
	\end{remark}
	
	\section{Discrete subgroups with critical exponent equal to 2}\label{section-delta2}
	
	In this section we give a proof of Theorem \ref{thm-delta2}, using the techniques developped in \cite{connellNaturalFlowCritical2023}. It uses the function $f$ defined in Lemma \ref{lem-f-unif-spsh}. We also outline a second proof, which involves the complete vector field $\mathfrak{X}=\nabla{f}$, called the \textit{natural flow} in \cite{connellNaturalFlowCritical2023}.
	
	\begin{proof}[Proof of Theorem \ref{thm-delta2}]
	Let $\Gamma$ be a discrete and torsion-free subgroup of $\PU(n,1)$ with critical exponent $\delta=2$ and assume that $\HNC/\Gamma$ contains a compact subvariety $A$ of positive dimension. 
	First, we remark that $\Gamma$ is non-elementary, as a consequence of Proposition \ref{prop-criteres}-(a) and (c). Thus $\Gamma$ admits a Patterson-Sullivan measure $(\mu_x)_{x\in \HNC}$. Fix a point $o\in \HNC$ and, for all $\theta\in \PHNC$, denote by $B_\theta:=B_{\theta}(\cdot,o)$ the Busemann function at $\theta$ which vanishes at $o$. Let ${f}$ be the $\Gamma$-invariant function defined on $\HNC$ by $f(x):=-\ln \lVert \mu_x\rVert$. By Lemma \ref{lem-f-unif-spsh}, for all tangent vector $v$ at a point $x\in \HNC$, we have
	\begin{equation*}
		i\partial\overline{\partial}{f}(v,Jv) \ge 0.
	\end{equation*}
		
	Let $\widetilde{A}\subset \HNC$ be the lift of $A$, $x$ be a regular point of $\widetilde{A}$ and $v$ be a non-zero vector in $T_x\widetilde{A}$. Then the plurisubharmonic function ${f}\vert_{ \widetilde{A}}$ is constant and consequently $i\partial\overline{\partial}{f}(v,Jv) = 0$.
	The inequality given by Lemma \ref{lem-f-unif-spsh} is thus an equality for this vector $v$. Using that Patterson-Sullivan measures are supported on $\Lambda(\Gamma)$, one sees that this equality can only happen when\begin{equation*}
		\forall \theta\in \Lambda(\Gamma), dB_{\theta}(v)^{2}+dB_{\theta}(Jv)^{2} = \lVert v\rVert^{2}.
	\end{equation*}
	This is possible if and only if $v\in \C v_{x\theta}$ for all $\theta\in \Lambda(\Gamma)$ and all $v\in T_x\widetilde{A}$, where $v_{x\theta} \in T_x\HNC$ is the unit vector at $x$ pointing in the direction of $\theta$. We deduce that $A$ has dimension 1 and that
	\begin{equation*}
		\forall \theta\in \Lambda(\Gamma),\ v_{x\theta}\in T_x\widetilde{A}.
	\end{equation*}
	Let $D$ be the unique complex geodesic containing $x$ for which $T_xD=T_x\widetilde{A}$. Then $\Lambda(\Gamma)\subset \partial D$, and hence the convex hull of $\Lambda(\Gamma)$ is contained in $D$. We deduce that $\Gamma$ preserves $D$. Moreover $A\subset D/\Gamma$ by Proposition \ref{prop-KH-hol-cvx}. To conclude, notice that $D/\Gamma$ is a Riemann surface containing a compact subvariety of positive dimension, so $D/\Gamma$ is compact and $\Gamma$ is a complex Fuchsian group. 
\end{proof}

\begin{remark}
		We now outline a second proof of Theorem \ref{thm-delta2}, inspired by \cite[Theorem 1.5]{connellNaturalFlowCritical2023}. According to \cite[Lemma 2.2]{connellNaturalFlowCritical2023}, the vector field $\mathfrak{X}$ defined by $\nabla{f}$ on $\HNC/\Gamma$ is complete, and defines a smooth flow $(\phi_t)_{t\in \R}$. Fix $x\in \HNC$ and let $\gamma$ be an integral curve for the vector field $\mathfrak{X}$. If $(Y_1,\dots,Y_k)$ spans a $k$-dimensional subspace $V$ of $T_x\HNC$, then one can define for all $t$ a $k$-dimensional subspace $V(t)$ of $T_{\gamma(t)}\HNC$ as the space spanned by $(\phi_{t\ast}Y_1,\dots,\phi_{t\ast}Y_k)$. Then the infinitesimal contraction rate of $V(t)$ is given by the real number $\tr(\nabla d{f}(x)\vert_{ V})$, see \cite[Lemma 2.5]{connellNaturalFlowCritical2023}.  For $x\in \HNC$ and $\theta\in \PHNC$, there is a real basis $(e_1,e_2,\dots,e_{2n})$ of $T_{x}\HNC$ with $e_1=v_{x\theta}$ the unit vector pointing in the direction of $\theta$ and $e_2=Je_1$, in which the matrix of $L_{\theta} + JL_{\theta}J^{-1}$ is
		$\operatorname{Diag}(2-\delta,2-\delta,2,\dots,2)$, see \cite[proof of Lemma 2.6]{connellNaturalFlowCritical2023}.
		In particular, if $\delta(\Gamma)=2$ and $V\subset T_x\HNC$ is a complex subspace, then\begin{equation*}
			\tr(L_{\theta\vert V}) = \frac{1}{2}\tr\left((L_{\theta} + JL_{\theta}J^{-1})\vert_{ V}\right) \ge 0,
		\end{equation*}
		with equality if and only if $V$ has complex dimension 1 and $V=\C v_{x\theta}$.
		Let $\widetilde{A}$ be the lift in $\HNC$ of a compact subvariety of positive dimension $A\subset\HNC/\Gamma$. Then we have for all regular point $x$ of $\widetilde{A}$\begin{equation*}
			\tr(\nabla d{f}\vert_{ T_x\widetilde{A}}) \ge \int_{\PHNC}\tr(L_{\theta\vert T_x\widetilde{A}})d\mu_x(\theta) \ge 0.
		\end{equation*}
		If, for some regular point $x$ of $\widetilde{A}$, the above inequality was strict, then $\phi_{-t}$ would contract $A$ for sufficiently small $t>0$, which would contradict the fact that $A$ is a volume minimizer in its homology class. Fixing from now on a regular point $x$ of $\widetilde{A}$, we deduce that $A$ has dimension 1 and, since Patterson-Sullivan measures of $\Gamma$ are supported on the limit set $\Lambda(\Gamma)$ of $\Gamma$, we get\begin{equation*}
			\forall \theta \in \Lambda(\Gamma), v_{x\theta} \in T_x\widetilde{A}.
		\end{equation*}
	We conclude as in the first proof.
\end{remark}
\vspace{-0.5cm}
\bibliographystyle{alpha}
\bibliography{biblio}
	
	\vspace{.1cm}
	\begin{small}\begin{tabular}{l}
			William, Sarem\\
			Univ. Grenoble Alpes, CNRS, IF, 38000 Grenoble, France\\
			william.sarem at univ-grenoble-alpes.fr
		\end{tabular}
	\end{small}		
\end{document}